\RequirePackage{fix-cm}
\documentclass[smallextended,envcountsect]{svjour3}

\usepackage[T1]{fontenc}
\usepackage[utf8]{inputenc}
\smartqed  
\usepackage{graphicx}
\usepackage{tikz-cd}
\usepackage{subfigure}
\usepackage{xcolor}

\usepackage{amsmath}
\usepackage{amssymb}
\usepackage{amscd}

\usepackage{url}
\usepackage{nicefrac}
\usepackage{mathrsfs}
\usepackage[outercaption]{sidecap}
\usepackage{amsfonts,epigraph}%
\usepackage{booktabs} 
\usepackage[all,cmtip]{xy}
\usepackage{enumerate}
\usepackage[misc]{ifsym}

\numberwithin{equation}{section}
\newtheorem{ass}[theorem]{Assumption}

\newtheorem{choices-notations}[theorem]{Choices and Notations}
\newtheorem{notation}[theorem]{Notation}

\newcommand{\exendproof}{\renewcommand{\qed}{\relax}\end{proof}}

\DeclareRobustCommand*{\nicefrac}[2]{\ifmmode\mathnicefrac{#1}
	{ #2}%
	\else\textnicefrac{#1}{#2}\fi}
\newcommand*{\textnicefrac}[2]{\check@mathfonts%
	\mbox{\raisebox{.5ex}{\fontsize\sf@size\z@\selectfont#1}\kern-.
		1em%
		/\kern-.1em\raisebox{- .25ex}{\fontsize\sf@size\z@\selectfont#2} }}
\newcommand*{\mathnicefrac}[2]{%
	\mathchoice
	{\m@fr@c{\scriptstyle}{#1}{#2}}
	{\m@fr@c{\scriptstyle}{#1}{#2}}
	{\m@fr@c{\scriptscriptstyle}{#1}{#2}}
	{\m@fr@c{\scriptscriptstyle}{#1}{#2}}}

\newcommand{\abs}[1]{\lvert#1\rvert}

\def\Ci{C^\infty}

\def\fequal#1{\stackrel{#1}{=}}

\def\lla{\langle}

\newcommand{\norm}[1]{\lVert#1\rVert}
\newcommand{\scalar}[2]{\lla#1,#2\rra}

\def\rra{\rangle}
\def\sqm1{\sqrt{-1}}

\def\tand{\mbox{\ \rm  and }}

\def\too{\longrightarrow}

\def\wt{\widetilde}
\def\x{\times}

\def\={\cong}
\def\>{\supset}
\def\<{\subset}

\def\12{\frac{1}{2}}
\def\0{^{\circ}}

\newcommand{\auindex}[1]{}

\def\CC{{\mathbb C}}
\def\EE{{\mathbb E}}

\def\NN{{\mathbb N}}

\def\RR{{\mathbb R}}

\def\Bb{{\mathcal B}}

\def\Dd{{\mathcal D}}

\def\a{\alpha}

\newcommand{\dop}{\operatorname{d}}

\def\F{\Phi}

\def\R{\RR}

\def\Si{\Sigma}

\DeclareMathOperator{\ad}{ad}

\newcommand\Calderon{Cal\-der{\'o}n}

 \DeclareMathOperator{\dist}{dist}

\DeclareMathOperator{\grad}{grad}

\DeclareMathOperator{\Graph}{\mathfrak G}
 
\DeclareMathOperator{\Id}{Id} \DeclareMathOperator{\image}{im}

 \DeclareMathOperator{\mmax}{max} \DeclareMathOperator{\mmin}{min}

\DeclareMathOperator{\ort}{ort}

\def\Vol{\mbox{\rm vol}}

\begin{document}

\title{Continuity of family of Calder{\'o}n projections
\thanks
{This work was supported by National Natural Science Foundation of China (Grant Nos. 11971245 and 11771331
), LPMC of MOE of China
and Nankai Zhide Foundation.}
}


\author{Bernhelm Boo{\ss}-Bavnbek\textsuperscript{1}
	\and Jian Deng\textsuperscript{2}
	\and Yuting Zhou\textsuperscript{3}
	\and Chaofeng Zhu\textsuperscript{3}
}

\authorrunning{B. Boo{\ss}-Bavnbek, J. Deng, Y. Zhou, and C. Zhu, \today} 

%

\institute{
  Chaofeng Zhu (\Letter) \email{zhucf@nankai.edu.cn}\\
  Bernhelm Boo{\ss}-Bavnbek \url{https://orcid.org/0000-0002-8865-7298} \email{booss@ruc.dk}\\
  Jian Deng \email{dengjian@cufe.edu.cn}\\
  Yuting Zhou \email{nkzhouyt@mail.nankai.edu.cn}\\
\at
{1} Department of Sciences and Environment, Roskilde University, DK-4000 Roskilde, Denmark
\at
{2} China Economics and Management Academy (CEMA),
           Central University of Finance and Economics,
           Beijing 100085, P. R. China
\at
{3} Chern Institute of Mathematics and LPMC,
           Nankai University,
           Tianjin 300071, P. R. China\\
}


\maketitle

\begin{abstract}
We consider a continuous family of linear elliptic differential operators of arbitrary order over a smooth compact manifold with boundary.
Assuming constant dimension of the spaces of inner solutions, we prove that the orthogonalized \Calderon\ projections of the underlying family of elliptic operators form a continuous family of projections. Hence, its images (the Cauchy data spaces) form a continuous family of closed subspaces in the relevant Sobolev spaces.
We use only elementary tools and classical results: basic manipulations of operator graphs and other closed subspaces in Banach spaces; elliptic regularity; Green's formula and trace theorems for Sobolev spaces; well-posed boundary conditions; duality of spaces and operators in Hilbert space; and the interpolation theorem for operators in Sobolev spaces.
\end{abstract}


\section{Introduction}\label{s:intro}
This paper provides a new approach to the investigation of \textit{Cauchy data spaces} (made of the normal traces at the boundary up to the order $d-1$ of the kernel of elliptic differential operators of order $d\ge 1$, that can be obtained as images of certain pseudo-differential projections over the boundary, called \textit{\Calderon\ projections}) under continuous or smooth variation of the underlying operators over a fixed manifold with boundary. The concept of the \textit{\Calderon\ projection} originated from
\textsc{\Calderon}'s observation in \cite{Cal63}.

Previous approaches to the variational problem were based either on purely functional-analytic, symplectic and topological arguments or on geometric and holomorphic analysis. For the first type of approach we refer to  \cite{BoFu98} that dealt with symmetric operators admitting self-adjoint Fredholm extensions and a certain unique continuation property (UCP). The variation was restricted to compact perturbations. By those assumptions, the authors achieved the continuous variation of the Cauchy data spaces in symplectic quotient Hilbert spaces, namely as Lagrangian subspaces.

The second type of approach is based on investigating spectral projec\-tions and exploiting the pseudo-differential calculus, via canonical and explicit constructions of Poisson operators and the \Calderon\ projection. See, e.g., the classical \cite{Ni95,Ni97} for Dirac type operators, based on the invertible double via gluing of \cite{BoWo93}, or our \cite{BoLe:2009,BoLeZh08,BCLZ} for arbitrary elliptic differential operators of first order with UCP, based on the ideas of general invertible doubles via a system of boundary value problems in \cite{Himpel-Kirk-Lesch:2004}.

Our present approach is a hybrid, changing repeatedly between the calculus of closed subspaces of the graph-theoretical approach and the geometric analysis of the \Calderon\ projections of the pseudo-differential approach. In that way we obtain the wanted generalization to linear elliptic differential operators of order $d\ge 1$ with weakened UCP requirements and, as a bonus, a much shorter path to the quoted results.

\subsection{Structure of the paper}
This paper consists of three sections. In this Section \ref{s:intro}, we explain the structure of the paper and state our main result.

In Section \ref{s:common-knowledge}, we fix the notations. The main topics are Sobolev spaces and domains of elliptic differential operators on manifolds with boundary; Green's forms; Cauchy data spaces; the homogenized Cauchy trace operator; the classical properties of the \Calderon\ projection; and \textsc{Neubauer}'s classical $\cap$ and $+$ \textit{arithmetic} of pairs of families of closed subspaces in Banach space \cite{Ne68}.

The proof of Theorem \ref{t:main} is in Section \ref{s:proof}. In Section \ref{ss:proof-for-s-ge-d-half}, assuming $s\ge \frac{d}{2}$, we obtain first the continuous variation of the solution spaces in the Sobolev space of order $s +\tfrac d2 $. By the continuity and surjectivity of the adjusted trace operator, that yields a continuous variation of the Cauchy data spaces in the Sobolev space of order $s$ over the boundary, and, furthermore, the continuous variation of  the family of $L^2$-orthogonalized
Calder{\'o}n projections in the operator norm of these Sobolev spaces.
This part of our results has been announced in \cite[Proof of Proposition 4.5.2, first part]{BoZh14}.
In Section \ref{ss:s<halfd}, we use the results of Section \ref{ss:proof-for-s-ge-d-half} to prove our theorem for $s<\frac{d}{2}$  by duality and interpolation property of spaces and operators in Sobolev scales. In the following Appendix, we show that the assumption about the constant dimension of the spaces of inner solutions
in Theorem \ref{t:main} can be weakened a little by finer analysis above.

\subsection{Our main result}

\begin{notation}\label{n:basic-notations}
Let $B$ be a topological space and $\mathscr{M}$ a compact smooth Riemannian manifold with boundary $\Si$. Let $\left(A_b\right)_{b\in B}$ be a family of linear elliptic differential operators of order $d\ge 1$, acting between sections of complex finite-dimensional Hermitian vector bundles $E,F$ over $\mathscr{M}$.

Let $\rho^d$ denote the {\em Cauchy trace operator}, mapping sections of $E$ over $\mathscr{M}$ to $d$-tuples of jets over $\Si$ in normal direction (these jets can be {\em adjusted}, i.e., homogenized to sections of the bundle $E'^d:=(E|_\Si)^d$, for details see Proposition \ref{p:trace}). Let
$$Z_{+,0}(A_b)\ :=\ \{u \in H^d(\mathscr{M}; E) \mid A_bu = 0 \tand \rho^d u=0\}$$
denote the space of all {\em inner solutions}. It is the finite-dimensional kernel of the {\em closed minimal realization} associated with $A_{b}$. Correspondingly, $Z_{-,0}(A_b) := Z_{+,0}(A^t_b)$
denotes the kernel of the closed minimal realization associated with the formal adjoint $A^t_{b}$.
\end{notation}


For the interesting case of Cauchy data spaces and $L^2$-orthogonalized (and so uniquely determined) \Calderon\ projections (see Section \ref{ss:weak-traces}), we shall prove

\begin{theorem}[Main result]\label{t:main}
Assume that
	\begin{enumerate}[(i)]\label{e:ucp-assumption}
	\item  for $s\ge \frac{d}{2}$\/, the two families of bounded extensions
	\[
	\bigl(A_{b, s+\frac{d}{2}}\colon
	 H^{s+\frac{d}{2}}(\mathscr{M};E) \too H^{s-\frac{d}{2}}(\mathscr{M};F)\bigr)_{b\in B}
	 \]
and
\[\bigl(A_{b, s+\frac{d}{2}}^t\colon H^{s+\frac{d}{2}}(\mathscr{M};F) \too H^{s-\frac{d}{2}}(\mathscr{M};E)\bigr)_{b\in B}
\]
	are continuous in the respective operator norms $\norm{\cdot}_{s+\frac{d}{2},s-\frac{d}{2}}$\/, and that the family of adjusted Green's forms (of Equation \eqref{e:J-adjusted}) $\bigl(\tilde J^t_{b,s} \colon H^s(\Si;F'^d) \to H^s(\Si;E'^d)\bigr)_{b\in B}$
	is continuous in the operator norm $\norm{\cdot}_{s,s}$\/;
	\item $\dim Z_{+,0}(A_b) \tand \dim Z_{-,0}(A_b)$  do not depend on $b\in B$.		
	\end{enumerate}
	Then for any $s\in \RR$, the family of $L^2$-orthogonalized Calder{\'o}n projections $\bigl(C^{\ort}_s(A_b)\bigr)_{b\in B}$ is continuous in the operator norm of the corresponding Sobolev space $H^s(\Si;E'^d)$.
\end{theorem}

\begin{remark}\label{r:continuity-assumptions} 
(a) Assumption (i) can be weakened by demanding continuous variation only
for $s\geq \frac{d}{2}$ and $s+\frac{d}{2}\in\NN$.
\par
(b) Let $M$ be covered by coordinate charts $(U,\varphi)\in \mathscr{A}$ for an atlas $\mathscr{A}$ with local trivializations $E|_{U},F|_{U}$. If for any $(U,\varphi)\in \mathscr{A}$, all multiple-order partial derivatives of the coefficients of the operators $A_b$ are continuous and uniformly bounded on $U\times B$ (cf. \cite[Section 1]{Atiyah-Singer:1971}), then Assumption (i) follows. 
\par
(c) In the literature on families  of elliptic operators over manifolds with boundary, strict weak inner unique continuation property is commonly assumed (that is, $Z_{\pm,0}(A_b)=\{0\}$).  Relaxing that assumption to  Assumption (ii), i.e., the constant dimensions of the spaces of inner solutions, was suggested in \cite{Himpel-Kirk-Lesch:2004}.
To prove the continuity of $\bigl(\ker A_{b,s+\frac{d}{2}}\bigr)_{b\in B}$ for $s\geq \frac{d}{2}$ (see Proposition \ref{p:kernel-cont-for-s-ge-dhalf}), we assume that the spaces $Z_{-,0}(A_b)$ are of finite constant dimension. Actually, the two statements are equivalent by Lemma \ref{l:closed-continuous}a. Once we have obtained the continuity of $\bigl(\ker A_{b,s+\frac{d}{2}}\bigr)_{b\in B}$, the assumption that the spaces $Z_{+,0}(A_b)$ are of finite constant dimension is equivalent to our conclusion
that the family of the images of the corresponding \Calderon\ projections is continuous (see Proposition \ref{p:Cauchy-traces-varying}).

In some special example, Assumption (ii) can be weakened. Please see the Appendix.
\par
(d) In \eqref{e:cauchy-data-space-basic} we define the \textit{Cauchy data space} $\Lambda_{-\frac{d}{2}}(A_b)  \< H^{-\frac{d}{2}}(\Si;E'^d)$ as the space of the homogenized Cauchy traces of the weak solutions $u$ of $A_bu=0$ and also in \eqref{e:cauchy-data-space-strong} the Cauchy data spaces $\Lambda_s(A_b)$
for $s\ge\frac{d}{2}$. According to Theorem \ref{t:Se69} and Corollary \ref{c:generalization-of-calderon-in-seeley}, these spaces are precisely the images of the corresponding \Calderon\ projections.

Clearly, the continuity of a family of projections of a Banach space in operator norm implies the continuity of their images in the gap topology, see also \cite[Section I.4.6]{Ka95}. Hence one can read Theorem \ref{t:main} as the claim of a continuous variation of the Cauchy data spaces depending on the parameter $b$ for each of these Sobolev orders $s$ --- under the assumption of constant dimensions of the spaces of inner solutions.
\par
(e) One of the most fundamental examples is the continuous variation of the Riemannian metric on a fixed smooth manifold, i.e., that in local coordinates all derivatives of the component functions of the metric vary continuously. Then the induced Laplace operators vary continuously in the sense of (b). The weak inner unique continuation property holds for Laplace operators. So both Assumptions (i) and (ii) are satisfied, and as a consequence the corresponding Cauchy data spaces vary continuously.
\par
(f) For applications of our results we refer to the spectral flow formulae for operator families with varying maximal domains as in \cite{BoZh14};
and, consequently, to the possibility of determining the precise number of negative eigenvalues in stability analysis of an essentially positive differential operator $A$ (appearing in the descriptions of, e.g., reaction-diffusion, wave propagation and scattering systems) by calculating the spectral flow of $\bigl((1-b)A+bA_+\bigr)_{b\in [0,1]}$, where $A_+$ is a suitably chosen strictly positive differential operator,  or more advanced expressions, cf. \cite{Latushkin-et-al:2018,BoZhu:2004,Zhu:2006} in the tradition of \textsc{Bott}'s Sturm type theorems \cite{Bo56}.
\end{remark}

\section{Main tools and notations of elliptic operators on manifolds with boundary} \label{s:common-knowledge}

Before 
proving the theorem, we fix the notations and recall the most basic concepts and tools.
We begin with a single operator.

\subsection{Our data}\label{ss:our-data}
\begin{enumerate}
	\item $\mathscr{M}$ is a smooth compact Riemannian manifold of dimension $n$ with boundary $\partial \mathscr{M}=:\Si$.
	\item $E,F\to \mathscr{M}$ are Hermitian vector bundles of fiber dimension $m$ with metric connections $\nabla^E$, $\nabla^F$. As in Notation \ref{n:basic-notations}, we set $E':=E|_\Si$ and $F':=F|_\Si$\/.
    \item $\Ci(\mathscr{M};E)$ denotes the space of smooth sections of $E$;
          $\mathscr{M}^\circ$ denotes the interior of $\mathscr{M}$, $\Ci_c(\mathscr{M}^\circ; E)$ denotes the space of smooth sections of $E$ with compact support in $\mathscr{M}^\circ$.
	\item $A\colon \Ci(\mathscr{M};E)\to\Ci(\mathscr{M};F)$ is an elliptic differential operator of order $d$.
	\item $A_0\colon \Ci_c(\mathscr{M}^\circ; E)\to\Ci_c(\mathscr{M}^\circ;F)$, where $A_0=A|_{\Ci_c(\mathscr{M}^\circ; E)}$.
	\item $A_0^t\colon \Ci_c(\mathscr{M}^\circ; F)\to\Ci_c(\mathscr{M}^\circ;E)$, where $A^t$ denotes the formal adjoint of $A$. 	
	\item  $A_{\mmin}:= \overline{A_0}$, $A_{\mmin}^t:= \overline{A_0^t}$, where we consider
$A_0\colon \Dd(A_0)\to L^2(\mathscr{M};F)$
as an unbounded densely defined operator from $L^2(\mathscr{M};E)$ to $L^2(\mathscr{M};F)$, and denote its closure by $\overline{A_0}$, see Section \ref{ss:sobolev}, in particular  Proposition \ref{p:minimal-domain} below. We write $A_{\mmax}:=(A_0^t)^{\ast}$\/, i.e.,
	\[
	\mathcal{D}(A_{\mmax})=\{u\in L^2(\mathscr{M};E)\mid Au\in L^2(\mathscr{M};F) \text{ in the distribution sense}\},
	\]
where $\mathcal{D}(\cdot)$ denotes the domain of an operator.
\end{enumerate}

Note that $A_{\mmin},A_{\mmax}$ are the closed \textit{minimal} and \textit{maximal extensions} of $A_0$. For a section $u\in \mathcal{D}(A_{\mmax})$,
the \textit{intermediate derivatives} $D^\a u$ (with $|\a|\le d$) need not exist as sections on $\mathscr{M}$, even though $Au$ does in the distribution sense, see \cite[Section 4.1, p. 61]{Grubb:2009}.

\subsection{The Sobolev scale and special relations for elliptic operators}\label{ss:sobolev}
For real $s$, we recall the definition of the Sobolev scale $H^s(\mathbf{M};\mathbf{E})$ for a complete smooth Riemannian manifold $\mathbf{M}$ without boundary. Then, for a compact manifold $\mathscr{M}$ with smooth boundary, the Sobolev scale is induced for non-negative $s$ by embedding and restriction.

We follow mostly \textsc{\Calderon} \cite[Section 3.1]{Cal76}, as reproduced and elaborated in \textsc{Frey} \cite[Chapters 0 and 1]{Frey2005On}, supplemented by \textsc{Lions} and \textsc{Magenes} \cite[Sections 1.7 and 1.9]{LM72} and \textsc{Tr{\`e}ves} \cite[Section III.2]{Treves:1}. We replace the regular subsets of $\RR^n$ in the classical literature by a smooth compact manifold with boundary embedded in a complete manifold without boundary. So, without restricting the general validity of our results, we assume, as we may, that
\begin{itemize}
	\item our compact Riemannian manifold $(\mathscr{M},g)$ \textit{with} boundary is embedded in a (metrically) complete smooth Riemannian manifold $( \mathbf{M},  \mathbf g)$  of the same dimension $n$ \textit{without} boundary,
	\item our bundles $E,F$ are extended to smooth Hermitian vector bundles $ \mathbf E, \mathbf F$ over $ \mathbf{M}$,
\item the elliptic differential operator $A$ is defined on $\mathscr{M}\cup \mathscr{N}$ where $\mathscr{N}$ denotes a collar neighbourhood of $\Si$ in $\mathbf M\setminus \mathscr{M}^{\circ}$.
\end{itemize}
%

\paragraph{Sobolev scale on complete manifolds without boundary.} First we recall the concept of the Sobolev scale for functions. The immediate generalization for sections of Hermitian bundles follows then.

On $ \mathbf{M}$ with Riemannian metric $ \mathbf g$, let $|\dop \Vol|$ denote the volume density derived from the metric.
Recall the \textit{Hodge--Laplace operator}
\begin{equation*}
\Delta_0^{ \mathbf{M}}\ :=\ \dop^t\dop\colon C_c^{\infty}( \mathbf{M})\too C_c^{\infty}( \mathbf{M}),
\end{equation*}
acting on functions, 
where $\dop^t$ denotes the formal adjoint of the exterior differential
$\dop\colon C^{\infty}( \mathbf{M})\to C^{\infty}( \mathbf{M};\Lambda^1( \mathbf M))$.
The operator $-\Delta_0^{ \mathbf{M}}$ is equal to the \textit{Laplace--Beltrami operator} on the Riemannian manifold $(\mathbf{M}, \mathbf g)$.
Let $L^2(\mathbf{M})$ denote the completion of $C_{c}^{\infty}( \mathbf{M})$ with respect to the norm induced by the $L^2$-inner product
\[
(u,v)_{L^2(\mathbf{M})}:=\int_{\mathbf{M}}u\bar{v}\,|\dop \Vol|,
\]
where $\bar{v}$ denotes the complex conjugate of $v$.
Since $\mathbf{M}$ is complete with respect to $\mathbf g$, $\Delta_0^{\mathbf{M}}$ is essentially self-adjoint (see \cite[Theorem 3]{Cordes:1972} or \cite[Section 3(A)]{Chernoff:1973}). So the closure of $\Delta_0^{ \mathbf{M}}$, $\Delta^{ \mathbf{M}}$, is a non-negative self-adjoint operator.
It gives rise to the \textit{Sobolev spaces} on $ \mathbf{M}$
\begin{equation}\label{e:sobolev-functions}
H^s( \mathbf{M})\ :=\ \mathcal{D}((\Delta^{ \mathbf{M}})^{s/2}),\ \ s\geq0,
\end{equation}
equipped with the graph norm.
By \cite[Theorem 1.1.2]{LM72}, we regain,
for $ \mathbf M=\RR^n$ and $s\in\NN\cup\{0\}$
the usual Hilbert space
\begin{equation*}
H^s( \mathbf{M})\ =\ \{u\in L^2( \mathbf M)\mid D^\alpha u \in L^2( \mathbf M) \text{ for }
|\alpha|\le s\},
\end{equation*}
where the partial differentiation $D^\alpha$ with multi-index $\alpha$ is applied in the distribution sense and the scalar product and norm are defined by
\[
\langle u,v\rangle_s\ :=\ \sum_{|\alpha|\le s} (D^\alpha u,D^\alpha v)_{L^2(\mathbf{M})}\
\tand\ \norm{u}_s\ :=\ \sqrt{\scalar{u}{u}_s}\,.
\]

For $s>0$, we define the space $H^{-s}(\mathbf{M})$ of distributions to be the so-called $L^2$-\textit{dual} of $H^s( \mathbf{M})$, i.e.,
\begin{equation}\label{e:L2-anti-dual}
H^{-s}( \mathbf{M}):=\{u\in\mathscr{D}'( \mathbf{M}) \mid \exists_c\forall_{v\in H^s(\mathbf{M})} \abs{u(\bar{v})} = \abs{\langle v,u\rangle_{s,-s}} \leq c \|v\|_{H^s(\mathbf{M} )}\},
\end{equation}
here $\langle v,u\rangle_{s,-s}:=\overline{u(\bar{v})}$ with a distribution $u\in\mathscr{D}'( \mathbf{M})$ acting on a test function $v$. Hence for
$u\in L^2(\mathbf{M})$ we have $\langle v,u\rangle_{s,-s}=
(v,u)_{L^2(\mathbf{M})}$, as nicely explained in \cite[Section 8.2]{Grubb:2009} and \cite[Section 1.1]{Gilkey:1995}.

The above constructions can be generalized for sections of any bundle $ \mathbf E\rightarrow  \mathbf{M}$ carrying an Hermitian structure $\mathbf{M}\ni p\mapsto \langle.,.\rangle|_{ \mathbf E_p}$ and an Hermitian connection.
Let
\begin{align*}
\nabla^{ \mathbf E}&\colon C^{\infty}( \mathbf{M}; \mathbf E)\longrightarrow C^{\infty}( \mathbf{M};T^* \mathbf{M}\otimes  \mathbf E)\ \tand\\
\nabla^{ \mathbf F}&\colon C^{\infty}( \mathbf{M}; \mathbf F)\longrightarrow C^{\infty}( \mathbf{M};T^* \mathbf{M}\otimes  \mathbf F)
\end{align*}
be Hermitian connections, i.e., connections that are compatible with the Hermitian metrics on $ \mathbf E$ and $ \mathbf F$ respectively. To define Sobolev spaces of sections in vector bundles, one replaces the Laplacian $\dop^t\dop$ in the previous definition \eqref{e:sobolev-functions} by the \textit{Bochner--Laplacian}s $(\nabla^{ \mathbf E})^t\nabla^{ \mathbf E}$ and $(\nabla^{ \mathbf F})^t\nabla^{ \mathbf F}$.

\paragraph{Sobolev scale on compact smooth manifolds with boundary.} For functions, the corresponding Sobolev space on the compact submanifold $\mathscr{M}$ with boundary $\Si$ is defined as the quotient
\begin{equation*}
H^s(\mathscr{M})\ :=\  H^s( \mathbf{M})/\left\{u\in H^s( \mathbf{M})\big|\  u|_{\mathscr{M}}=0\right\}, s\in\R\tand s\ge 0.
\end{equation*}
In other words, $H^s(\mathscr{M})$ coincides algebraically with the space of restrictions to $\mathscr{M}^{\circ}$ of the elements of $H^s(\mathbf{M})$. The norm of $H^s(\mathscr{M})$ is given by the quotient norm, that is,
\[\norm{u}_{H^s(\mathscr{M})}=\inf\norm{U}_{H^s(\mathbf{M})} \ \ \text{for all $U\in H^s(\mathbf{M})$ with $U=u$ a.e. on $\mathscr{M}^{\circ}$.}
\]
In our smooth case, the definition coincides with the interpolation $H^s(\mathscr{M}) = [H^m(\mathscr{M}) , H^0(\mathscr{M})]_{\theta}$,
$(1-\theta)m=s$, $m$ integer, $0\leq\theta\leq1$.
See \cite[Theorems 1.9.1 and 1.9.2]{LM72}. For $s\ge 0$, an important subspace is the function space
$H_0^s(\mathscr{M}):=\overline{C_c^{\infty}(\mathscr{M}^{\circ})}^{\|.\|_{H^s(\mathscr{M})}}$.

More generally and quite similarly, we can  define Sobolev spaces of sections in vector bundles such as  $H^s(\mathscr{M};E)$ and
\begin{equation}\label{e:H_0}
H_0^s(\mathscr{M};E)\ :=\ \overline{C_c^{\infty}(\mathscr{M}^{\circ};E|_{\mathscr{M}^{\circ}})}^{\|\cdot\|_{H^s(\mathscr{M};E)}} \ \text{for $s\in\RR$, $s\ge 0$}.
\end{equation}

\paragraph{Sobolev scale on closed manifolds.}
For any Hermitian vector bundle $G$ over the closed manifold $\Sigma$,
we can define the Sobolev spaces $H^s(\Sigma;G)$ for all $s\in \RR$ as in \cite[Section 8.2]{Grubb:2009} or \cite[Section 1.3]{Gilkey:1995}. 
Note that $C^{\infty}(\Sigma;G)$ is dense in $H^s(\Sigma;G)$ for all $s\in\RR$.
Then the $L^2$-scalar product for smooth sections can be extended to a \textit{perfect pairing} between $H^s(\Sigma;G)$ and $H^{-s}(\Sigma;G)$ for all $s\in\RR$.
That is, from \cite[Lemma 1.3.5(e)]{Gilkey:1995}, the pairing $(f,h)_{L^2(\Si;G)}$ extends continuously to a perfect pairing
\[
H^s(\Sigma;G)\times H^{-s}(\Sigma;G)\too \CC,
\]
which we denoted by $\langle \cdot,\cdot\rangle_{s,-s}$ in \eqref{e:L2-anti-dual}.
\begin{remark}\label{r:perfect-pairing}
	Let $H, K$ be Hilbert spaces. A bounded sesquilinear form $\F\colon H\times K\to\CC$ is called a \textit{perfect pairing} if it induces on each of $H, K$ an isomorphism to the dual of the other. More precisely, we obtain
the induced conjugate linear map from $K$ to the space of bounded linear functionals on $H$ by
\[
v\mapsto \F(\cdot,v)\ \ \text{for $v\in K$,}
\]
the induced linear map from $H$ to the space of bounded conjugate linear functionals on $K$ by
\[
u\mapsto \F(u,\cdot)\ \ \text{for $u\in H$.}
\]
Both are isomorphisms; moreover, the functionals are bounded by
\begin{equation*}
\norm{v}_K\ =\ \sup_{0\ne u\in H}\frac{\abs{\F(u,v)}}{\norm{u}_H}\ \tand\
\norm{u}_H\ =\ \sup_{0\ne v\in K}\frac{\abs{\F(u,v)}}{\norm{v}_K}\,.
\end{equation*}
\end{remark}


\begin{notation}
Let $X,Y$ be normed spaces, we denote the normed algebra of bounded linear operators from $X$ to $Y$ by $\Bb(X,Y)$; for $X=Y$, $\Bb(X):=\Bb(X,X)$.
We use shorthand $\norm{\cdot}_s$ for the norm in $H^s(\cdot;\cdot)$, $s\in\RR$; and $\norm{\cdot}_{r,s}$ for the operator norm in $\Bb(H^r(\cdot;\cdot),H^s(\cdot;\cdot))$, $r,s\in\RR$.
\end{notation}

\paragraph{Special relations for elliptic operators.}
We fix the notation, in particular the sign conventions.
Let
$d_g(\cdot,\cdot)$ be the \textit{distance function}; (locally) it is the arc length of the minimizing geodesic.
In a collar neighbourhood of
$\Si$ in $\mathscr{M}$, say $V$, the function
\[
V\ni p\mapsto x_1(p)\ :=\ d_g(p,\Sigma),\   p \in \mathscr{M}
\]
is smooth and defines the \textit{inward unit normal field} $\nu:=\grad x_1$ and \textit{inward unit co-normal field} $\nu^\flat:=\dop x_1$.

Let $T^*\mathscr{M}$ denote the \textit{cotangent vector bundle} of $\mathscr{M}$, $S(\mathscr{M})$ the \textit{unit sphere bundle} in $T^*\mathscr{M}$ (relative to the Riemannian metric $g$), and $\pi \colon S(\mathscr{M})\rightarrow  \mathscr{M}$ the projection. Then associated with any linear differential operator $A$ of order $d$ there is a vector bundle homomorphism
$$\sigma_d(A)\colon \pi^*E\rightarrow \pi^*F\/,$$
which is called the \textit{principal symbol} of $A$. In terms of local coordinates, $\sigma_d(A)$ is obtained from $A$ by replacing $\partial/\partial x_j$ by $\mathrm{i}\xi_j$ in the highest order terms of $A$ (here $\xi_j$ is the $j$th coordinate in the cotangent bundle). $A$ \textit{elliptic} means that $\sigma_d( A)$ is an isomorphism.

For elliptic operators there is an important relation (the \textit{G{\aa}rding inequality}) between the graph norm, originating from the basic $L^2$ Hilbert space, and the corresponding Sobolev norm. More precisely, we recall from \cite[Proposition 1.1.1]{Frey2005On}

\begin{proposition}\label{p:minimal-domain}
	Assume that $A$ is an elliptic operator of order $d$. Then
	\begin{enumerate}[(a)]
		\item The graph norm of $A$ restricted to $\Ci_c(\mathscr{M}^\circ;E)$ is equivalent to the Sobolev norm $\norm{\cdot}_{H^d(\mathscr{M};E)}$\/.
		\item In particular, $\mathcal{D}(A_{\mmin})=H^d_0(\mathscr{M};E)$ and $\mathcal{D}(A^t_{\mmin})=H^d_0(\mathscr{M};F)$.
		\item $H^d(\mathscr{M};E)\< \mathcal{D}(A_{\mmax})$ is dense.
	\end{enumerate}
\end{proposition}

\subsection{Green's formula, traces of Sobolev spaces over the boundary, and weak traces for elliptic operators }\label{ss:traces}
Let $j\in \NN\cup\{0\}$.
Let $\gamma^j \colon C^{\infty}(\mathscr{M};E)\rightarrow C^{\infty}(\Sigma;E')$ denote the trace map $\gamma^j u:= (\nabla_{\nu}^E)^ju|_{\Sigma}$ yielding the $j$th jet in normal direction. Set
\begin{equation}\label{e:rho-d}
\rho^d\ :=\ \left(\gamma^0,...,\gamma^{d-1}\right)\colon C^{\infty}(\mathscr{M};E)\too C^{\infty}({\Si};E'^d).
\end{equation}
Analogously, $\nabla^F$ gives rise to trace maps $\gamma^j \colon C^{\infty}(\mathscr{M};F)\rightarrow C^{\infty}(\Sigma;F')$. The corresponding maps for $F$ will also be denoted by $\gamma^j$ and $\rho^d$.

We recall \textit{Green's Formula}, e.g., from \textsc{Seeley} \cite[Equation 7]{See66}, \textsc{Tr\`eves} \cite[Equation III.5.41]{Treves:1}, \textsc{Grubb} \cite[Proposition 11.3]{Grubb:2009}, or \textsc{Frey} \cite[Proposition 1.1.2]{Frey2005On}, with a description of the operator $J$ in the error term:
\begin{proposition}[Green's Formula for differential operators of order $d\ge 1$]\label{Green's formula}
Let $A \colon C^{\infty}(\mathscr{M};E)\too C^{\infty}(\mathscr{M};F)$ be a linear differential operator of order $d$. Then 	
	there exists a (uniquely determined) differential operator
	\begin{equation*}
	J \colon C^{\infty}(\Sigma;E'^d)\too C^{\infty}(\Sigma;F'^d),
	\end{equation*}
	such that for all $u\in C^{\infty}(\mathscr{M};E), v\in C^{\infty}(\mathscr{M};F)$ we have
	\begin{equation}\label{e:Green}
	( Au,v)_{L^2(\mathscr{M};F)}- (u,A^tv)_{L^2(\mathscr{M};E)}\ =\ (J\rho^du,\rho^dv)_{L^2(\Sigma;F'^d)}.
	\end{equation}
	$J$ is a matrix of differential operators $J_{kj}$ of order $d-1-k-j$, $0\leq k,j\leq d-1$, and $J_{kj}=0$ if $k+j>d-1$ ($J$ is upper skew-triangular). Moreover, for $j=d-1-k$ we have explicitly given homomorphisms
	\begin{equation}\label{e:greens-form-explicit}
	J_{k,d-1-k}\ =\ \mathrm{i}^d (-1)^{d-1-k}\sigma_d(A)(\nu^\flat).
	\end{equation}
\end{proposition}

\begin{remark}\label{r:greens-form}
(a)	For $d=1,2,3$, we visualize the structure of the matrix $J$,
	\[
	\begin{pmatrix} J^{[0]}_{00}\end{pmatrix},\
	\begin{pmatrix} J^{[1]}_{00}&  J^{[0]}_{01}\\
	J^{[0]}_{10} &   0\end{pmatrix},\
	\begin{pmatrix} J^{[2]}_{00}&  J^{[1]}_{01}&  J^{[0]}_{02}\\
	J^{[1]}_{10} & J^{[0]}_{11} &   0\\
	J^{[0]}_{20} & 0 &   0\end{pmatrix},\ \text{etc.},
	\]
	where the orders of the differential operators of the entries were marked by a superscript $[\langle order\rangle]$.
\newline
(b) From the explicit form of the skew diagonal elements of
$J$ in \eqref{e:greens-form-explicit}, we get that $J$
is invertible for any elliptic operator $A$.
\newline
(c) If $J$ is invertible, Green's Formula \eqref{e:Green} extends to $(u,v)\in \mathcal{D}(A_{\mmax})\times H^d(\mathscr{M};F)$, where the right-hand side is interpreted as the $L^2$-dual pairing
$$\oplus_{j=0}^{d-1}H^{-d+j+\frac{1}{2}}(\Sigma;F')\ \times\ \oplus_{j=0}^{d-1}H^{d-j-\frac{1}{2}}(\Sigma;F')\ \too\ \CC.$$
\end{remark}
%
%

With \cite[Theorem 1.1.4]{Frey2005On}, we obtain a slight reformulation, sharpening, and generalization of the classical \textit{Sobolev Trace Theorem} (see also \cite[Section 9.1]{Grubb:2009} and \cite[Lemma 16.1]{Tar07}):

\begin{proposition}[Sobolev Trace Theorem]\label{p:trace}
	\begin{enumerate}[(1)]
	\item We have continuous trace maps $\rho^d$ (obtained by continuous extension):
	\begin{eqnarray*}
		(a) &&\rho^d\colon H^{d+s}(\mathscr{M};E)\too \displaystyle\oplus_{j=0}^{d-1}H^{d+s-j-\frac{1}{2}}(\Sigma;E') \text{ for $s>-\frac{1}{2}$}\, ,\\
		(b) &&\rho^d\colon\mathcal{D}(A_{\mmax})\too \displaystyle\oplus_{j=0}^{d-1}H^{-j-\frac{1}{2}}(\Sigma;E').
	\end{eqnarray*}
	Moreover, the map (a) is surjective and has a continuous right-inverse $\eta^d$.
	\item If $u\in\mathcal{D}(A_{\mmax})$, then $u\in H^d(\mathscr{M};E)$ if and only if
	$$\rho^du\ \in\ H^{d-\frac{1}{2}}(\Sigma;E')\oplus\cdots\oplus H^{\frac{1}{2}}(\Sigma;E').$$
	\item For $\rho^d$ on $\mathcal{D}(A_{\mmax})$ we have
	$
	\ker\/\rho^d\ = \ H_0^d(\mathscr{M};E).
	$	
\end{enumerate}
\end{proposition}

\begin{remark}\label{r:cauchy-trace-map}
(a) 	Following \textsc{Grubb} \cite[Section 9.1]{Grubb:2009}, we call the preceding scale of operators $\rho^d$ with domains in different Sobolev spaces by one name: the \textit{Cauchy trace operator} associated with the order $d$.
\newline
(b)	It is well known that the trace operators do not extend to the whole $L^2(\mathscr{M};E)$. For the special case of the half-space in $\RR^n$, it is shown in \cite[Remark 9.4]{Grubb:2009} that the 0-trace map $\gamma^0$ makes sense on $H^s(\RR^n_+)$ if and only if $s > \12$. The Cauchy trace operator $\rho^d$ extends, however, to $\mathcal{D}(A_{\mmax})$, though in a  way that depends on the choice of the elliptic operator $A$ but with $\ker\rho^d=H^d_0(\mathscr{M};E)$ independent of $A$.  That makes the claims 1b and 2 of the preceding proposition particularly interesting.
\newline	
(c) Claim 3 admits replacing the abstract definition of $H_0^d(\mathscr{M};E)$ in \eqref{e:H_0} by
a concrete check of the Cauchy boundary data of a given section.
The inclusion $H_0^d(\mathscr{M};E)\<\ker\bigl(\rho^d|_{\mathcal{D}(A_{\mmax})}\bigr)$ is obvious. There remains to be shown that the converse holds. Via local maps and claim 2, this reduces to the Euclidean case (see \cite[Theorem 9.6]{Grubb:2009}).
\end{remark}

\paragraph{Homogenization of Sobolev orders.}
Following \textsc{\Calderon}\ \cite[Section 4.1, p. 76]{Cal76} and using the notation of \cite[p. 26]{Frey2005On}, we introduce a \textit{homogenized} (\textit{adjusted}) Cauchy trace operator $\wt\rho^d$.
We set $\Delta^{E'}:=(\nabla^{E'})^t\/\nabla^{E'}$, where $\nabla^{E'}$ denotes the
restriction of $\nabla^{E}$ 
to $\Sigma$. Since $\Delta^{E'}+1$ is a positive symmetric elliptic differential operator of second order on the closed manifold $\Sigma$, it possesses a discrete spectral resolution (e.g., \cite[Section 1.6]{Gilkey:1995}). Then $\Phi :=(\Delta^{E'}+1)^{1/2}$ is a pseudo-differential operator of order 1 which induces an isomorphism of Hilbert spaces
\[
\Phi_{(s)} \colon H^s(\Si;E')\ \too\ H^{s-1}(\Si;E') \text{ for all $s\in\RR$}
\]
and, in fact, generates the Sobolev scale $H^s(\Sigma;E')$, see the Sobolev scale of an unbounded operator in \cite[Section 2.A]{BrLe01}.

In order to achieve that all boundary data are of the same Sobolev order, we introduce the matrix
\begin{equation*}
\Phi_d\quad :=\quad\left(
\begin{array}{cccc}
\Phi^{\frac{d-1}{2}} & 0 & \cdots & 0 \\
0 & \Phi^{\frac{d-3}{2}} & \cdots & 0 \\
\vdots & \vdots & \ddots & \vdots \\
0 & 0 & \cdots & \Phi^{\frac{-d+1}{2}}\\
\end{array}
\right).
\end{equation*}

We set $\wt{\rho}^d:=\Phi_d\circ\rho^d$, $\wt{\eta}^d:=\eta^d\circ\Phi_d^{-1}$. So, we obtain a condensed and adjusted Trace Theorem as a corollary to Proposition \ref{p:trace}:

\begin{corollary}[Homogenized trace map]\label{c:tilde-rho} 	
We have continuous trace maps $\wt\rho^d$ (obtained by continuous extension):
\begin{enumerate}[(a)]
	\item $\wt\rho^d\colon H^{s+\frac{d}{2}}(\mathscr{M};E)\ \too\ H^s(\Sigma;E'^d)$ for $s>\frac{d}{2}-\frac{1}{2}$,
	\item $\wt\rho^d\colon \mathcal{D}(A_{\max})\ \too\ H^{-\frac{d}{2}}(\Sigma;E'^d)$.
\end{enumerate}
Furthermore, the map (a) is surjective and has a continuous right-inverse, $\wt\eta^d$. \end{corollary}

\begin{remark}
(a)	We have
	$\mathcal{D}(A_{\min})= \ker\bigl(\wt\rho^d\colon \mathcal{D}(A_{\max})\to H^{-\frac{d}{2}}(\Sigma;E'^d)\bigr).$
\newline
(b) Analogous constructions for the bundle $F$ lead to $\Delta^{F'}\colon C^{\infty}(\Sigma;F')\to C^{\infty}(\Sigma;F')$. Whenever this causes no ambiguity we will denote the corresponding matrices,
$\Phi^{F'}$ and $\Phi_d^{F'}$, again by $\Phi$ and $\Phi_d$, respectively.
\end{remark}
%

We can replace the boundary operator $J$ of Green's Formula (Proposition \ref{Green's formula}) by its adjusted version
\begin{equation}\label{e:J-adjusted}
\wt{J}\ :=\ (\Phi_d^{F'})^{-1}\circ J\circ (\Phi_d^{E'})^{-1}.	
\end{equation}

It follows that all components of $\wt{J}$,
\[
\wt{J}_{ij}=\Phi^{\frac{2i+1-d}{2}}J_{ij}\Phi^{\frac{2j+1-d}{2}}
\]
are pseudo-differential operators of order $i+j+(1-d)+(d-1)-i-j=0$.
So for any $s\in\RR$, \begin{equation*}
\wt J\colon\ H^s(\Si,E'^d)\ \too\ H^s(\Si,F'^d).
\end{equation*}
$\wt J$ is upper skew-triangular, with invertible elements on the skew diagonal for elliptic operator $A$.

For later use we give the adjusted version of Green's Formula. 
It is valid for arbitrary linear differential operators of order $d\ge 1$ with smooth coefficients acting between sections of smooth Hermitian vector bundles $E$ and $F$ over a smooth compact Riemannian manifold $\mathscr{M}$ with boundary $\Si$:
With the preceding notations of $\wt J$ for the adjusted Green boundary operator of \eqref{e:J-adjusted} and $\wt\rho^d$ for the adjusted Cauchy trace operators of Corollary \ref{c:tilde-rho}a,  we have
\begin{equation}\label{e:green-adjusted}
\ (Au,v)_{L^2(\mathscr{M};F)}-(u,A^tv)_{L^2(\mathscr{M};E)}\ =\ (\wt J\wt \rho^du,\wt \rho^dv)_{L^2(\Sigma;F'^d)}
\end{equation}
for $s\geq\frac{d}{2}$ and all $u\in H^{s+\frac{d}{2}}(\mathscr{M};E)$, $v\in H^{s+\frac{d}{2}}(\mathscr{M};F)$. For elliptic differential operators, \eqref{e:green-adjusted} remains valid 
for $u\in \mathcal{D}(A_{\mmax})$, $v\in H^d(\mathscr{M};F)$ in the following extended version
\begin{equation*}
 (Au,v)_{L^2(\mathscr{M};F)}- (u,A^tv)_{L^2(\mathscr{M};E)}\ =\ \lla\wt J\wt \rho^du,\wt \rho^dv\rra_{-\frac{d}{2},\frac{d}{2}}\,,
\end{equation*}
where $\wt \rho^du\in H^{-\frac{d}{2}}(\Si;E'^d)$ and $\wt \rho^dv\in H^{\frac{d}{2}}(\Si;F'^d)$.
	
\subsection{Properties of the \Calderon\ projection}\label{ss:weak-traces}

We give our version of perhaps not widely known classical results concerning the \Calderon\ projection.
%
%

In his famous note \cite{Cal63} of 1963, \textsc{A. \Calderon} introduced the concept of a pseudo-differential projection onto the Cauchy data of solutions of a system of homogeneous elliptic differential equations over smooth compact manifolds with boundary, later called the \Calderon\ projection, see below Theorem \ref{t:Se69} and Corollary \ref{c:generalization-of-calderon-in-seeley}.
While the \textit{Cauchy trace operator} was introduced for \textit{sections} in Equation \eqref{e:rho-d} and Corollary \ref{c:tilde-rho}, we define the \textit{Cauchy data spaces} for elliptic differential operators.
We recall

\begin{definition}[Cauchy Data Spaces]\label{d:cauchy-data-spaces}
	Let $A, d, \mathscr{M}, \Si, E, F$ be given as in Section \ref{ss:our-data}.
Based on the homogenized Cauchy trace operators in Corollary \ref{c:tilde-rho}, we define the Cauchy data spaces as follows:
\newline
(a)
\begin{equation}\label{e:cauchy-data-space-basic}
\Lambda_{-\frac{d}{2}}(A)\ :=\ \{h\in H^{-\frac{d}{2}}(\Sigma;E'^d)\mid  \exists u \in \ker A_{\max} \text{ with }\wt{\rho}^du=h\},
\end{equation}		
i.e., the space of boundary values of weak solutions to $A$ (= sections belonging to the maximal domain of $A$ that vanish under the operation of $A$ in the distribution sense). 
\newline
(b) 
For $s\ge \frac{d}{2}$\/, 
\begin{equation}\label{e:cauchy-data-space-strong}
\Lambda_{s}(A):=\{h\in H^s(\Sigma;E'^d)\mid  \exists u \in H^{s+\frac{d}{2}}(\mathscr{M};E),
Au=0 \text{ with }\wt{\rho}^du=h\}.
\end{equation}
\end{definition}
\begin{remark}
(a) For later use, we rewrite the Cauchy data spaces:
\begin{gather}
   \Lambda_{-\frac{d}{2}}(A)\ = \ \wt{\rho}^d(\ker A_{\mmax}), \quad  \tand \notag\\
  \Lambda_{s}(A)= \wt{\rho}^d(\ker A_{s+\frac{d}{2}})= \wt\rho^d\bigl(\ker A_{\mmax}\cap H^{s+\frac{d}{2}}(\mathscr{M};E)\bigr)\ \text{ for $s\ge \frac{d}{2}$\/,} \label{e:cauchy-data-spaces-s-ge=d/2}
\end{gather}
where
$A_{s+\frac{d}{2}} \colon H^{s+\frac{d}{2}}(\mathscr{M};E)\to H^{s-\frac{d}{2}}(\mathscr{M};F)$ is the extension of $A$.
In fact,
for $u\in  H^{s+\frac{d}{2}}(\mathscr{M};E)$, $v\in \Ci_c(\mathscr{M}^\circ;F)$,
we have
\[\bigl(A_{s+\frac{d}{2}}u,v\bigr)_{L^2(\mathscr{M};F)}\ =\ (u,A^t_0v)_{L^2(\mathscr{M};E)},\] so
$A_{\mmax}|_{H^{s+\frac{d}{2}}(\mathscr{M};E)}=A_{s+\frac{d}{2}}$, thus \eqref{e:cauchy-data-spaces-s-ge=d/2} follows.
\newline
(b)
For a subspace $V$ of $L^2(\Si;F'^d)$, we denote by \[V^{\bot {L^2}}\ :=\ \{g\in L^2(\Si;F'^d)\mid (f,g)_{L^2(\Si;F'^d)}=0 \ \text{for all $f\in V$}\}.\]
Let $\wt{J}$ be the adjusted Green's form of \eqref{e:J-adjusted}.
Now we recall an important relationship of the Cauchy data spaces between $A$ and $A^t$:
\begin{equation}\label{e:cauchy-data-space-transposed}
  \Lambda_{\frac{d}{2}}(A^t)\ =\ \left(\wt{J}\Lambda_{\frac{d}{2}}(A)\right)^{\bot {L^2}}\cap\,  H^{\frac{d}{2}}(\Si;F'^d).
\end{equation}
It was proved in \cite[Proposition 2.1.1]{Frey2005On} and will be used in the proof of our Corollary
\ref{c:l2-orthogonality-of-calderon-in-seeley}.
\end{remark}

There is a vast literature on pseudo-differential operators and the symbolic calculus. We shall only draw on the general knowledge regarding pseudo-differential operators over \textit{closed} manifolds.
Let $k\in \RR$. Let $\Psi_k(\Si;G_1,G_2)$ denote the space of all $k$th order pseudo-differential operators mapping sections of a smooth vector bundle $G_1$
to sections of a smooth vector bundle $G_2$ over the same closed manifold $\Si$; for $P\in \Psi_k(\Si;G_1,G_2)$, let $\sigma_k(P)$ denote the principle symbol of $P$.

\begin{proposition}\label{p:pseudo-differential-property}(cf. \cite[Lemmas 1.34 and 1.35]{Gilkey:1995})\\
(a) For $P\in \Psi_k(\Si;G_1,G_2)$, there is a unique formal adjoint of $P$, denoted by $P^t$, such that
\[
(Pf,h)_{L^2(\Sigma;G_2)}=(f,P^th)_{L^2(\Sigma;G_1)}\ \ \text {for all $f\in C^{\infty}(\Sigma;G_1),h\in C^{\infty}(\Sigma;G_2)$},
\]
and $P^t\in \Psi_k(\Si;G_2,G_1)$, $\sigma_k(P^t)=\sigma_k(P)^*$.\\
(b) If $Q\in \Psi_k(\Si;G_1,G_2), P\in \Psi_e(\Si;G_2,G_3)$, then $PQ\in \Psi_{k+e}(\Si;G_1,G_3)$ and
$\sigma_{k+e}(PQ)=\sigma_k(P)\sigma_e(Q)$.\\
(c) Continuity property with respect to Sobolev spaces:
Each $P\in \Psi_k(\Si;G_1,G_2)$ extends to a continuous linear map form
$H^{s+k}(\Si;G_1)$ to $H^s(\Si;G_2)$ for all real $s$.
\end{proposition}

With the preceding notations we recall the classical knowledge about \Calderon\ projections.

\begin{theorem}[A. Calder\'{o}n 1963; R. T. Seeley 1966, 1969] \label{t:Se69}
	Let $A$ be an elliptic differential operator of order $d$ over a smooth compact Riemannian manifold $\mathscr{M}$ with boundary $\Si$, acting between sections of Hermitian vector bundles $E,F$ over $\mathscr{M}$.
	Then there exists a (in general, not uniquely determined) zeroth order classical pseudo-differential operator, called the {\em Calder\'{o}n projection} of $A$, $C(A)=C_\infty(A)\colon \Ci(\Sigma;E'^d)\to C^{\infty}(\Sigma;E'^d)$, such that
	$C(A)$ is idempotent, i.e., $C(A)^2=C(A)$ and
	\begin{equation}\label{e:range-of-calderon}
	\image C_{-\frac{d}{2}}(A) \ = \ \wt\rho^d(\ker A_{\mmax})\/.	
\end{equation}
	Here, for every $s\in\RR$, we denote the extended projection by $C_s(A)\colon H^s(\Sigma;E'^d)$ $\to H^s(\Sigma;E'^d)$\/. 
\end{theorem}

\begin{remark}\label{r:calderon-sources}
	A detailed construction of the Calder\'{o}n projection and a careful proof of \eqref{e:range-of-calderon} 
can be found in \cite[Section 2.3]{Frey2005On}, originally from
	\textsc{Seeley}'s \cite{Palais-Seeley:1965,See66,Seeley:1968}  and \textsc{Calder\'{o}n}'s \cite{Cal76} and inspired by \textsc{H{\"o}rmander}'s \cite{Ho66}; for the properties see also \cite[Section 11.1]{Grubb:2009} and, in the special case $d=1$, \cite[Chapters 12-13]{BoWo93} and \cite[Section 5]{BoLeZh08} .
\end{remark}

\begin{remark}\label{r:regular-wellposed-boundary}
	The Calder\'{o}n projections are very useful in the treatment of boundary value problems for elliptic differential operators. We recall some relevant definitions and properties directly from \cite{Frey2005On}.
	\par
	(a) Assume $P\in \Psi_0(\Si;E'^d,E'^d)$ idempotent. To consider $P$ as a
boundary condition we associate with it the  \textit{realisation}
$A_P\colon \mathcal{D}(A_{P})\to L^2(\mathscr{M};F)$ with
	\[
	\mathcal{D}(A_{P}):=\{u\in  H^d(\mathscr{M};E)\mid P\wt{\rho}^du=0\}.
	\]
	We define the \textit{weak domain} by
	\[
	\mathcal{D}(A_{\mmax,P}):=\{u\in \mathcal{D}(A_{\mmax})\mid P\wt{\rho}^du=0\}.
	\]
	\par
	(b) Following \cite[Definition 1.2.5]{Frey2005On},
	we call $P$  a \textit{regular} boundary condition if
	\[
	\mathcal{D}(A_P)=\mathcal{D}(A_{\mmax,P}).
	\]
	$P$ is called \textit{well-posed} if it is regular and $\image A_P$ has finite codimension.

 From \cite[Proposition 2.1.2]{Frey2005On} we recall equivalent conditions:
	\begin{itemize}\label{i:boundary-condition-Fredholm}
		\item $P$ is a regular boundary condition if and only if $A_P \colon \mathcal{D}(A_P)\to L^2(\mathscr{M};F)$ is left-Fredholm, i.e., $\dim\ker A_P< \infty$ and $\image A_P$ closed.
		\item The boundary condition $P$ is well-posed if and only if $A_P\colon\mathcal{D}(A_P)\to L^2(\mathscr{M};F)$ is Fredholm.
	\end{itemize}
  \par(c) With these notations we obtain from \textsc{Seeley}'s achievements in \cite{Seeley:1968}, reproduced and worked out in \cite[Theorem 2.1.4(ii)]{Frey2005On} the following \textit{operational conditions} on the principle symbols of boundary pseudo-differential operators:
Regularity and moreover well-posedness hold respectively if and only if for all $q\in \Sigma$, $\xi\in T^*_q\Sigma$, $\xi\neq 0$
	\[
\sigma_0(P)(q,\xi)\colon  \image \sigma_0\bigl(C(A)\bigr)(q,\xi)\too E'^d|_{q} \text{ injective},
	\]
respectively,
	\[
	\sigma_0(P)(q,\xi)\colon  \image \sigma_0\bigl(C(A)\bigr)(q,\xi)\too \image \sigma_0(P)(q,\xi)  \text{ invertible.}
	\]
In particular,
since $\sigma_0(\Id)=\Id$ and
	\[
	\sigma_0\bigl(C(A)\bigr)(q,\xi)\colon  \image \sigma_0\bigl(C(A)\bigr)(q,\xi) \too \image \sigma_0\bigl(C(A)\bigr)(q,\xi)
	\]
	is just the identity,
$\Id$ is a regular boundary condition and $C(A)$ is a well-posed boundary condition.
\par
	(d)
	For regular $P$ we have a lifting jack for regularity, also called \textit{higher regularity} (see \cite[Theorems 2.2.1 and 2.2.3]{Frey2005On}):
Let $s\in \NN \cup\{0\}$.
Assume that $u\in L^2(\mathscr{M};E)$ satisfies
\[Au\in H^s(\mathscr{M};F),\ \ \ \ P\wt \rho^du \in H^{s+\frac{d}{2}}(\Sigma;E'^d).\]
Then $u\in H^{s+d}(\mathscr{M};E)$.
When $P$ is also well-posed, then this regularity argument holds for all real $s\geq0$.
\par
(e) For later use we give the following simple description of $\image C_s(A)$ for all
$s\in \RR$\, :
\begin{equation}\label{e:calderon2cauchy-data-spaces-any}
 \image C_s(A)\  
 =\ \image C_t(A)\, \cap\, H^s(\Si;E'^d)\ \text{ for any real $t\leq s$.}
\end{equation}
Since $s\ge t$, $H^s(\Si;E'^d)\subset H^t(\Si;E'^d)$, then $\image C_s(A) =
 \image\left(C_t(A)|_{H^s(\Si;E'^d)}\right)$.
Hence, the inclusion $\<$ of (\ref{e:calderon2cauchy-data-spaces-any}) is trivial.
For the opposite inclusion we exploit that $C(A)$ is a zeroth order pseudo-differential \textit{idempotent}: so, for any $C_t(A)f\in  H^s(\Si;E'^d)$, we have
\[C_t(A)f\ =\ C^2_t(A)f\ =\ C_s(A)C_t(A)f.\]
In particular, for all
$s\ge -\frac{d}{2}$\,:
\begin{equation}\label{e:calderon2cauchy-data-spaces}
 \image C_s(A)\
 =\ \image C_{-\frac{d}{2}}(A)\, \cap\, H^s(\Si;E'^d)\ =\ \Lambda_{-\frac{d}{2}}(A)\, \cap\, H^s(\Si;E'^d).
\end{equation}
As a side result, we obtain for all $s\in\RR$
\begin{equation}\label{e:chain-of-cds}
\image C_s(A)\ =\ \overline{\image C_\infty(A)}^{\|\cdot\|_{H^s(\Sigma;E'^d)}}=\ \overline{\image C_t(A)}^{\|\cdot\|_{H^s(\Sigma;E'^d)}}\ \text{for any real $t\geq s$}.
\end{equation}
\end{remark}

According to the preceding Remark \ref{r:regular-wellposed-boundary}e, the image of $C_s(A)$,
$s\in\RR$, does not depend on the choices of \Calderon\ projection $C(A)$ in the \Calderon --Seeley Theorem \ref{t:Se69}. Moreover we can prove the following generalization of Equation \eqref{e:range-of-calderon}.
\begin{corollary}\label{c:generalization-of-calderon-in-seeley}
The validity of the claim \eqref{e:range-of-calderon} in Theorem \ref{t:Se69} holds also for Sobolev orders $s\ge \frac{d}{2}$, yielding
	\begin{equation}\label{e:range-of-calderon-all}
	\image(C_s(A)) = \wt\rho^d\bigl(\ker A_{\mmax}\cap H^{s+\frac{d}{2}}(\mathscr{M};E)\bigr)\ \text{ for $s\ge \frac{d}{2}$}.
	\end{equation}
\end{corollary}
\begin{proof}
 By the Sobolev Trace Theorem (cf. Corollary \ref{c:tilde-rho}) and Equations \eqref{e:range-of-calderon} and \eqref{e:calderon2cauchy-data-spaces}, we obtain
	\[
	\image C_s(A)\ \>\ \wt\rho^d\bigl(\ker A_{\mmax}\cap H^{s+\frac{d}{2}}(\mathscr{M};E)\bigr) \text{ for $s> \frac{d}{2}-\frac{1}{2}$}\,,
	\]
	i.e., the inclusion $\>$ of \eqref{e:range-of-calderon-all}, actually for a wider range of $s$ than claimed.
	\par
	Now we turn to the proof of the inclusion $\<$  for $s\geq \frac{d}{2}$.
 First by the preceding Remark \ref{r:regular-wellposed-boundary}c, $C(A)$ is a well-posed boundary condition. Let $s\geq \frac{d}{2}$.
	If $f\in \image C_s(A)\fequal{\eqref{e:calderon2cauchy-data-spaces}}\Lambda_{-\frac{d}{2}}(A)\, \cap\, H^s(\Si;E'^d)$, then there is a $u\in \ker A_{\mmax}$, such that $f=\wt\rho^d u$, so $C(A)\wt\rho^d u=\wt\rho^d u \in H^s(\Si;E'^d)$. By the higher regularity for well-posed boundary conditions of the preceding Remark \ref{r:regular-wellposed-boundary}d,
	we have $u\in H^{s+\frac{d}{2}}(\mathscr{M};E)$, so $f\in \wt\rho^d\bigl(\ker A_{\mmax}\cap H^{s+\frac{d}{2}}(\mathscr{M};E)\bigr)$.
	Thus we get \eqref{e:range-of-calderon-all}.\qed
\end{proof}
\begin{remark}
In Definition \ref{d:cauchy-data-spaces}, we defined the Cauchy data spaces for $s=-\frac{d}{2}$ and $s\ge \frac{d}{2}$\/.
By Theorem \ref{t:Se69} and Corollary \ref{c:generalization-of-calderon-in-seeley}, these spaces coincide with the images of the extensions $C_s(A)$ of the Calder{\'o}n projection for those $s$. Moreover by \eqref{e:chain-of-cds} and Corollary \ref{c:generalization-of-calderon-in-seeley}, the images of $C_s(A)$ for all $s\in\RR$ are unique determined by the Cauchy data spaces of the elliptic operator.
That motivates us to define the Cauchy data spaces $\Lambda_s(A)$ for all $s\in\RR$ as the images of $C_s(A)$, yielding a \textit{chain of Cauchy data spaces}
\begin{equation}\label{e:cauchy-data-spaces}
  \Lambda_s(A)\ :=\ \image C_s(A) \ \text{for all $s\in\RR$}.
\end{equation}
\end{remark}

By \eqref{e:cauchy-data-space-transposed} and Corollary \ref{c:generalization-of-calderon-in-seeley},
we immediately have
\begin{equation}\label{e:range-of-calderon-transposed}
  \image C_{\frac{d}{2}}(A^t) \ =\ \left(\wt{J}\bigl(\image C_{\frac{d}{2}}(A)\bigr)\right)^{\bot {L^2}}\cap\,  H^{\frac{d}{2}}(\Si;F'^d).
\end{equation}

Now we can prove the following generalization of Equation \eqref{e:range-of-calderon-transposed}.
The preceding corollary and the following one will be used in Section \ref{s:proof} in the proof of our Main Theorem (Theorem \ref{t:main}).
\begin{corollary}\label{c:l2-orthogonality-of-calderon-in-seeley}

There is an $L^2$-orthogonal decompositions of complementary closed subspaces
	\begin{equation}\label{e:orthogonal-decomposition}
	\image C_s(A)\ \oplus^{\bot {L^2}}\ \wt{J}^t\left(\image C_s(A^t)\right) \ =\ H^s(\Si;E'^d) \text{ for $s\ge 0$},
	\end{equation}	
	where $\wt J$ is defined as in (\ref{e:J-adjusted}).
\end{corollary}

Our proof of Corollary \ref{c:l2-orthogonality-of-calderon-in-seeley} below
needs the (uniquely determined) \textit{$L^2$-orthogonalized \Calderon\ projection} which we are going to introduce now -- and use extensively later in our Section \ref{s:proof}.

\begin{lemma}\label{l:calderon-ort}
Let $C:=C(A)$ be a \Calderon\ projection as introduced in Theorem \ref{t:Se69}, i.e., a zeroth order classical pseudo-differential operator,
\[
C(A)\colon C^{\infty}(\Sigma;E'^d)\too C^{\infty}(\Sigma;E'^d)
\]
with $C^2=C$ and $\image C_{-\frac{d}{2}}(A)=\wt{\rho}^d(\ker A_{\mmax})$.
Then there exists a unique zeroth classical pseudo-differential operator $C^{\ort}=C^{\ort}(A)\colon C^{\infty}(\Sigma;E'^d)\to C^{\infty}(\Sigma;E'^d)$ with
\begin{equation}\label{e:orthogonal-calderonprojector-property}
(C^{\ort})^2=C^{\ort},\ \ CC^{\ort}=C^{\ort}, \ \ C^{\ort}C=C,
\end{equation}
with self-adjoint $L^2$-extension $C^{\ort}_0(A)$ on $H^0(\Si;E'^d)$ and with
\begin{equation}\label{e:invariant-L2orthogonal-Calderon}
\image C^{\ort}_s(A)= \image C_s(A)\ \ \ \text{for all $s\in\RR$}.
\end{equation}
\end{lemma}
\begin{proof}[of the lemma]
Since $CC^t+(\Id-C^t)(\Id-C)$ is a formally self-adjoint elliptic pseudo-differential operator with trivial kernel, it is invertible.
As in \cite[Lemma 12.8]{BoWo93}, we define $C^{\ort}$ by
\[
C^{\ort}:=CC^t\left(CC^t+(\Id-C^t)(\Id-C)\right)^{-1},
\]
and infer that it
is still a classical pseudo-differential operator of order $0$. Moreover, it is symmetric
\begin{equation}\label{e:orthogonal-calderonprojector-symmetry}
(C^{\ort}f,h)_{L^2(\Sigma;E'^d)}\ =\ (f,C^{\ort}h)_{L^2(\Sigma;E'^d)},\ \ \text{for all $f,h\in  C^{\infty}(\Sigma;E'^d)$}
\end{equation}
which implies (over the closed manifold $\Si$) that its $L^2$-extension is self-adjoint. The symmetry property \eqref{e:orthogonal-calderonprojector-symmetry} and
the algebraic equalities of \eqref{e:orthogonal-calderonprojector-property} follow as in loc. cit. by calculation, then the invariance of the range in \eqref{e:invariant-L2orthogonal-Calderon}
and the uniqueness of $C^{\ort}$ follow.\qed
\end{proof}

\begin{proof}[of Corollary \ref{c:l2-orthogonality-of-calderon-in-seeley}]
 For $s\geq\frac{d}{2}$\/, $\image C_s(A)\ni f=\wt{\rho}^du$ and $\image C_s(A^t)\ni g=\wt{\rho}^dv$ with $u\in\ker A_{\mmax}\cap H^{s+\frac{d}{2}}(\mathscr{M};E), v\in\ker A^t_{\mmax}\cap H^{s+\frac{d}{2}}(\mathscr{M};F)$, we obtain
\begin{eqnarray*}
  (f,\wt{J}^tg)_{L^2(\Si;E'^d)}\ &=&\ (\wt{J}\wt{\rho}^du,\wt{\rho}^dv)_{L^2(\Si;E'^d)}\\ &\fequal{\eqref{e:green-adjusted}}&
(Au,v)_{L^2(\mathscr{M};F)} - (u,A^tv)_{L^2(\mathscr{M};E)}=0-0=0.
\end{eqnarray*}
That proves the $L^2$-orthogonality in \eqref{e:orthogonal-decomposition}.
\par	
Now let $f\in H^s(\Si;E'^d)$ for $s=\frac{d}{2}$. 
To prove \[f \in \image C_{\frac{d}{2}}(A)+ \wt{J}^t\bigl(\image C_{\frac{d}{2}}(A^t)\bigr),\] i.e., the claimed decomposition in \eqref{e:orthogonal-decomposition} for $s=\frac{d}{2}$, we rewrite
\begin{equation}\label{e:orthogonal-decomposition-d/2}
f\ =\ C^{\ort}_{\frac{d}{2}}(A)f+f-C^{\ort}_{\frac{d}{2}}(A)f.
\end{equation}
We shall show that $f-C^{\ort}_{\frac{d}{2}}(A)f\in \wt{J}^t\bigl(\image C_{\frac{d}{2}}(A^t)\bigr)$. First we observe that
\begin{equation}\label{e:step1}
f-C^{\ort}_{\frac{d}{2}}(A)f\in \bigl(\image C_{\frac{d}{2}}(A)\bigr)^{\bot {L^2}}\cap H^{\frac{d}{2}}(\Si;E'^d).
\end{equation}
In fact, $h'\in\image C_{\frac{d}{2}}(A)=\image C^{\ort}_{\frac{d}{2}}(A)$ implies that $h'= C^{\ort}_{\frac{d}{2}}(A)h$ for some $h\in H^{\frac{d}{2}}(\Si;E'^d)$. Then we have
\begin{multline*}
\bigl(f-C^{\ort}_{\frac{d}{2}}(A)f,h'\bigr)_{L^2(\Si;E'^d)}\ =\
\bigl(f-C^{\ort}_{\frac{d}{2}}(A)f,C^{\ort}_{\frac{d}{2}}(A)h\bigr)_{L^2(\Si;E'^d)}\\
=\ \bigl(C^{\ort}_{\frac{d}{2}}(A)(f-C^{\ort}_{\frac{d}{2}}(A)f),h\bigr)_{L^2(\Si;E'^d)}\ =\ (0,h')_{L^2(\Si;E'^d)}\ =\ 0.\end{multline*}
Next from the fact that $\wt{J}$ is an invertible zeroth
	pseudo-differential operator, we obtain for any $h\in H^{\frac{d}{2}}(\Si;E'^d)$,
	\begin{eqnarray*}
		0 \ &\fequal{\eqref{e:step1}}&\  \bigl(f-C^{\ort}_{\frac{d}{2}}(A)f, C^{\ort}_{\frac{d}{2}}(A)h\bigr)_{L^2(\Si;E'^d)} \\
		&=&\ \bigl(f-C^{\ort}_{\frac{d}{2}}(A)f, \wt J^{-1} \wt J C^{\ort}_{\frac{d}{2}}(A)h\bigr)_{L^2(\Si;E'^d)} \\
		&=&\ \bigl((\wt J^t)^{-1}(f-C^{\ort}_{\frac{d}{2}}(A)f),  \wt J C^{\ort}_{\frac{d}{2}}(A)h\bigr)_{L^2(\Si;F'^d)},
	\end{eqnarray*}
	so
\begin{equation*}\label{e:step2-3}
(\wt J^t)^{-1}\bigl(f-C^{\ort}_{\frac{d}{2}}(A)f\bigr)\in \left(\wt{J}\bigl(\image C_{\frac{d}{2}}(A)\bigr)\right)^{\bot {L^2}}\!\cap  H^{\frac{d}{2}}(\Si;F'^d)
\fequal{\eqref{e:range-of-calderon-transposed}} \image C_{\frac{d}{2}}(A^t),
\end{equation*}
	and we are done for $s=\frac{d}{2}$.

	For $s>\frac{d}{2}$, the $L^2$-complement in $H^s(\Si;E'^d)$ of \eqref{e:orthogonal-decomposition} follows from the preceding result for $s=\frac{d}{2}$ and the facts that $C(A),C(A^t)$ and $C^{\ort}$ are pseudo-differential projections of order zero and $\wt{J}$ is an invertible
pseudo-differential operator of order zero.
	Finally,
	\eqref{e:orthogonal-decomposition} holds for $0\leq s< \frac{d}{2}$, since $H^{\frac{d}{2}}(\Si;E'^d)$ is dense in $H^s(\Si;E'^d)$ and $H^s(\Si;E'^d)\subset L^2(\Si;E'^d)$.
	\qed\end{proof}

\begin{remark}\label{r:symplectic-form}
(a) In the proof of Corollary \ref{c:l2-orthogonality-of-calderon-in-seeley}, we got that
\[\ker C^{\ort}_{s}(A)=\image \bigl(\Id -C^{\ort}_{s}(A)\bigr)=\wt{J}^t\bigl(\image C_{s}(A^t)\bigr)=\wt{J}^t\bigl(\Lambda_s(A^t)\bigr)\  \text{ for $s\ge 0$},\] which will be used in Section \ref{ss:proof-for-s-ge-d-half}.
\newline
(b) For a symmetric elliptic differential operator $A$ and $s\geq 0$, $\wt J^t$ defines a (strong) symplectic form on ($L^2(\Si;E'^d)$) $H^s(\Si;E'^d)$ with the Cauchy data space as a Lagrangian subspace according to \eqref{e:orthogonal-decomposition}.
\end{remark}

\subsection{\textsc{Neubauer}'s arithmetic of families of closed linear subspaces in Banach space}

We refer to a functional analysis fact from our \cite[Appendix A.3]{BoZh14}
regarding the continuity of families of closed subspaces in Banach space. We restate it in the following lemma that is based on \textsc{Neubauer}'s elementary, but deeply original \cite{Ne68}. We impose the gap topology on the space of closed linear subspaces of a given Banach space.
We recall the concept of the gap between subspaces and the quantity $\gamma$ ("angular distance") that is useful in our estimates.

\begin{definition}\label{d:gap}(cf. \cite[Sections IV.2.1 and IV.4.1]{Ka95})
Let $X$ be a Banach space.
\newline
(a) Denote by $S_M$ the unit sphere of $M$ for any closed linear subspace $M$ of $X$. For any two closed linear subspaces $M,N$ of $X$, we set
	\begin{eqnarray*}
		\delta(M,N)&:=&\left\{
		\begin{array}{ll}
			\sup_{u\in S_M}\dist(u,N), & \hbox{if $M\neq \{0\}$,} \\
			0, & \hbox{if $M=\{0\}$.}
		\end{array}
		\right.\\
		\hat{\delta}(M,N)&:=& \max\{\delta(M,N),\delta(N,M)\}.
	\end{eqnarray*}
	$\hat{\delta}(M,N)$ is called the \textit{gap} between $M$ and $N$.
\newline
We set
\begin{eqnarray}\label{e:mimimal-gap}
		\gamma(M,N)&:=&\left\{
		\begin{array}{ll}
			\inf_{u\in M \setminus N}\frac{\dist(u,N)}{\dist(u,M\cap N)}, & \hbox{if $M\not\subset N$,} \\
			1, & \hbox{if $M\subset N$.}
		\end{array}
		\right.
\end{eqnarray}
(b) We say, a sequence $(M_n)_{n=1,2\dots}$ of closed linear subspaces \textit{converges} to $M$ if $\hat{\delta}(M_n,M)\rightarrow 0$ for $n\rightarrow \infty$. We write $M_n\to M$. Correspondingly we declare when a mapping $M$ from a topological space $B$ to the space of closed subspaces shall be called \textit{continuous} at $b_0\in B$.
\end{definition}

\begin{remark}\label{r:gap}
Denote the set of all closed operators from $X$ to $Y$ by $\mathcal{C}(X,Y)$. If $A_1,A_2\in \mathcal{C}(X,Y)$, their graphs $\Graph(A_1):=\{(x,A_1x)\in X\times Y \mid x\in \Dd(A_1)\}$, $\Graph(A_2)$ are closed linear subspaces in the product Banach space $X\times Y$. We use the gap $\hat{\delta}(\Graph(A_1),\Graph(A_2))$ to measure the "distance" between $A_1$ and $A_2$.
Obviously, for $A',A\in \mathcal{B}(X,Y)$, we have (cf. \cite[Theorem IV.2.14]{Ka95})
\begin{equation*}
  \hat{\delta}(\Graph(A'),\Graph(A))\leq \|A'-A\|.
\end{equation*} 
\end{remark}
\begin{lemma}\label{l:gamma-closed-positive}(cf. \cite[Theorem IV.4.2]{Ka95})
Let X be a Banach space and let $M,N$ be closed subspaces of $X$.
In order that $M+N$ be closed, it is necessary and sufficient that $\gamma(M,N)$>0.
\end{lemma}

\begin{lemma}\label{l:closed-continuous}(cf. \cite[Poposition A.3.13 and Corollary A.3.14]{BoZh14})
	Let $X$ be a Banach space and let $\left(M_b\right)_{b\in B}, \left(N_b\right)_{b\in B}$ be two families of closed subspaces of $X$, where $B$ is a parameter space. Assume that $M_{b_0}+N_{b_0}$ is closed for some $b_0\in B$, and $(M_b)_{b\in B}$, $(N_b)_{b\in B}$ are both continuous at $b_0$\,.
\newline
(a) Then $(M_b\cap N_b)_{b\in B}$ is continuous at $b_0$ if and only if $(M_b+N_b)_{b\in B}$ is continuous at $b_0$.
	\newline
	(b) Assume furthermore that for $b\in B$,
	$\dim (M_b\cap N_b)\equiv$ constant $<+\infty$ or $\dim X/(M_b+N_b)\equiv$ constant $<+\infty$. Then the families $\left(M_b\cap N_b\right)_{b\in B}$ and $\left(M_b+N_b\right)_{b\in B}$ are both continuous at $b_0$.
\end{lemma}
\begin{remark}\label{r:closes-subspaces-sum}
For better understanding the proof of the preceding lemma in \cite{BoZh14}, note that:
according to \cite[Corollary A.3.12b]{BoZh14} and \cite[Theorem IV.4.2]{Ka95}), if $M_{b_0}+N_{b_0}$ is closed and $(M_b)_{b\in B}$, $(N_b)_{b\in B}$ and $(M_b\cap N_b)_{b\in B}$ are all continuous at $b_0$, then we get that $M_b+N_b$ is closed in a whole neighbourhood of $b_0$ in $B$.
\end{remark}
%
%
%
%

\section{Proof of our main theorem}\label{s:proof}
We shall divide the proof of Theorem \ref{t:main} into two cases, both under the assumption of
constant dimensions of the spaces of inner solutions $Z_{+,0}(A_b)= \ker A_{b,\mmin}$ and $Z_{-,0}(A_b)= \ker A^t_{b,\mmin}$\/:
\begin{enumerate}
  \item[(1)] In Section \ref{ss:proof-for-s-ge-d-half} we deal with the case $s\geq \frac{d}{2}$ in three steps. (i) In Proposition \ref{p:kernel-cont-for-s-ge-dhalf}, we obtain that $\bigl(\ker A_{b,s+\frac{d}{2}}\bigr)_{b\in B}$ is continuous in  $H^{s+\tfrac d2}(\mathscr{M};E)$. 
      (ii) Since the Cauchy trace map is bounded and surjective for $s\geq \frac{d}{2}$\/, we can deduce the main achievement of this Subsection, Proposition \ref{p:Cauchy-traces-varying}, and obtain that
the family $\bigl(\wt\rho^d(\ker A_{b,s+\frac{d}{2}})\bigr)_{b\in B}$ is continuous in $H^{s}(\Si;E'^d)$. That means $\bigl(\image C^{\ort}_s(A_b)\bigr)_{b\in B}$ is continuous in $H^{s}(\Si;E'^d)$. 
(iii) So, according to Corollary \ref{c:sufficient-condition-families-of-projections}, we can conclude that the corresponding family $\bigl(C^{\ort}_s(A_b)\colon H^s(\Si;E'^d)\hookleftarrow\bigr)_{b\in B}$ of \Calderon\ projections is continuous in the operator norm for all $s\geq \frac{d}{2}$. That proves our Main Theorem for such $s$.
 \item[(2)]  We use the results of case (1) (i.e., $s\geq \frac{d}{2}$) to show that for $s<\frac{d}{2}$ the family $\bigl(C^{\ort}_s(A_b)\colon H^s(\Si;E'^d)\hookleftarrow\bigr)_{b\in B}$ of \Calderon\ projections is continuous in the operator norm by duality and interpolation property of spaces and operators in Sobolev scales.
That is the content of Section \ref{ss:s<halfd}.
\end{enumerate}	
\begin{remark}\label{r:L2-orthogonalization-calderon}
We emphasize that all the \Calderon\ projections in this section are assumed to be $L^2$-orthogonalized, that is, with Lemma \ref{l:calderon-ort},
\[C=C^{\ort},\ \ \ \ \ C_s(A)=C^{\ort}_s(A) \ \text{for $s\in\RR$}.\]
\end{remark}

%
%
\subsection{Proof of our main theorem for $s\ge \frac{d}{2}$}
\label{ss:proof-for-s-ge-d-half}
We recall some of the technical ingredients and results obtained previously in our \cite[Proposition 4.5.2]{BoZh14}.

\begin{ass}\label{a:continuous-family-for-s-ge-dhalf}
Let $s\ge \frac{d}{2}$\/. We assume that the family

\[
\bigl(A_{b,s+\frac{d}{2}} \colon H^{s+\tfrac d2}(\mathscr{M};E) \to H^{s-\tfrac d2}(\mathscr{M};F)\bigr)_{b\in B}
\]
is a continuous family in the operator norm $\norm{\cdot}_{s+\frac{d}{2},s- \frac{d}{2}}$\/.
\end{ass}

For dealing with the case $s\ge \frac{d}{2}$, we introduce the following definition.
%
%
\begin{notation}
 Based on Remark \ref{r:regular-wellposed-boundary}a, for $s\geq\frac{d}{2}$, we denote by $A_{s+\frac{d}{2},P}$ the operators
\[
A_{s+\frac{d}{2},P}\colon\{u\in H^{s+\frac{d}{2}}(\mathscr{M};E)\mid P\wt \rho^du=0\}\too H^{s-\tfrac d2}(\mathscr{M};F),
\]
for any boundary condition $P\colon C^{\infty}(\Si;E'^d)\to C^{\infty}(\Si;E'^d)$.
We write shorthand $A_P:=A_{d,P}$\/.
\end{notation}

%

For any elliptic operator $A$ over a smooth compact manifold with boundary, recall $A_{\mmin}\colon H^d_0(\mathscr{M};E)\to L^2(\mathscr{M};F)$. It is well known and was emphasised above in Notation \ref{n:basic-notations} that $\ker A_{\mmin}$ consists only of smooth sections and is finite-dimensional. That follows from the interior regularity for elliptic operators (e.g., \cite[Theorem 5.11.1]{Taylor96}) and one can use the interior elliptic estimate to prove that $A_{\mmin}$ is left-Fredholm, i.e., $\dim \ker A_{\mmin}< +\infty$ and $\image A_{\mmin}$ is closed (e.g., \cite[Propositions 1.1.1 and A.1.4]{Frey2005On}). Later we shall use the following slight generalization:

\begin{lemma}\label{l:s-Amin-A-semifredholm} 
For $s\geq \frac{d}{2}$,
$\ker A_{s+\frac{d}{2},\Id}=\ker A_{\mmin}$ is finite-dimensional and consists of smooth sections.
\end{lemma}

\begin{proof}
We only need to prove the equality.
By Proposition \ref{p:trace}(3), $\Dd(A_{\Id})=H^d_0(\mathscr{M};E)$, so $A_{\Id}=A_{\mmin}$.
As just emphasized, we have
 $\ker A_{\mmin}\subset \{u\in C^{\infty}(\mathscr{M};E)\mid Au=0 \tand \wt\rho^du=0\}$.
Obviously we have for $s\geq \frac{d}{2}$
\[\{u\in C^{\infty}(\mathscr{M};E)\mid Au=0 \tand \wt\rho^du=0\}\subset\ker A_{s+\frac{d}{2},\Id}\subset\ker A_{\mmin}.\]
So we get the equality.
\qed
\end{proof}

In the following lemma we prove that $\image A_{s+\tfrac d2}$ is closed in $H^{s-\tfrac d2}(\mathscr{M};F)$ and get information about the quotient space
$H^{s-\tfrac d2}(\mathscr{M};F)/ \image A_{s+\tfrac d2}$ .

\begin{lemma}\label{l:strongdecomposition}
For $s\ge \frac{d}{2}$, there is an $L^2$-orthogonal decomposition of complementary closed subspaces
\begin{equation}\label{e:s>d/2-L^2-decomposition}
 H^{s-\tfrac d2}(\mathscr{M};F)\ =\ \image A_{s+\tfrac d2}\oplus^{\bot {L^2}}\ker A^t_{\mmin}.
\end{equation}
\end{lemma}
\begin{proof}
The $L^2$-orthogonality follows directly from Green's Formula (Proposition \ref{Green's formula}) and in adjusted form \eqref{e:green-adjusted}. In fact, for $s\ge\frac{d}{2}$,
$u\in H^{s+\tfrac d2}(\mathscr{M};E)$ and $v\in \ker A^t_{\mmin}$\/, we have
\[(Au,v)_{L^2(\mathscr{M};F)}\ =\ (u, A^tv)_{L^2(\mathscr{M};E)}+(\wt J \wt \rho^du,\wt \rho^dv)_{L^2(\Si;F'^d)}\ = 0.\]

Next we prove
\begin{equation}\label{e:L^2-decomposition}
  L^2(\mathscr{M};F)\ =\ \image A_C\oplus \ker A^t_{\mmin},
\end{equation}
where $C:=C^{\ort}(A)$ denotes the $L^2$-orthogonalized \Calderon\ projection defined in Lemma \ref{l:calderon-ort}.

By Remarks \ref{r:regular-wellposed-boundary}c, b, the \Calderon\ projection is a well-posed
boundary condition; hence $A_C \colon \mathcal{D}(A_C)\to L^2(\mathscr{M};F)$ is Fredholm, where
$\mathcal{D}(A_C)=\{u\in H^d(\mathscr{M};E) \mid C\wt{\rho}^d u=0\}$.
Thus $\image A_C$ is closed, then we have the decomposition $L^2(\mathscr{M};F) =\image A_C\oplus\ker\,(A_C)^*$, where we consider $A_C$ as an unbounded densely defined operator from $L^2(\mathscr{M};E)$ to $L^2(\mathscr{M};F)$ and denote its adjoint by $(A_C)^*$.
So \eqref{e:L^2-decomposition} will follow from
\begin{equation}\label{e:ker-calderon-boundary}
  \ker\,(A_C)^*\ =\ \ker A^t_{\mmin}\/.
\end{equation}

Now we shall prove \eqref{e:ker-calderon-boundary}.
In fact, according to \cite[Proposition 1.2.6]{Frey2005On},
\begin{equation*}
(A_C)^*= A^t_{\mmax,C^{\ad}} \text{ with }  C^{\ad}:=(\wt J^t)^{-1}(\Id-C^t)\wt J^t.
\end{equation*}
Note that $C^{\ad}\in \Psi_0(\Si;F'^d, F'^d)$ is idempotent and defines a well-posed boundary condition for $A^t$:
According to our assumption $C=C^{\ort}$, we have $C=C^t$. By Corollary \ref{c:l2-orthogonality-of-calderon-in-seeley}, we have
\[\wt J^t \bigl(\image C_{\frac{d}{2}}(A^t)\bigr)= \image \bigl(\Id-C^{\ort}_{\frac{d}{2}}(A)\bigr).\] Then we get
$\image C_{\frac{d}{2}}(A^t)=\image C_{\frac{d}{2}}^{\ad}$, where $C_{\frac{d}{2}}^{\ad}\colon H^{\frac{d}{2}}(\Si;F'^d)\to H^{\frac{d}{2}}(\Si;F'^d)$. Thus
\begin{equation}\label{e:calderon-adjoint-boundary-condition}
  C_{\frac{d}{2}}^{\ad} \colon \image C_{\frac{d}{2}}(A^t) \to \image C_{\frac{d}{2}}^{\ad}
\end{equation}
is just the identity.
So by \cite[Proposition 2.1.2]{Frey2005On}, $C^{\ad}$ is a well-posed boundary condition for $A^t$.
Then by Remark \ref{r:regular-wellposed-boundary}b, we have $A^t_{\mmax,C^{\ad}}=A^t_{C^{\ad}}$, thus
\begin{eqnarray*}
  \ker A^t_{\mmax,C^{\ad}}\ &=&\ \ker A^t_{C^{\ad}}  \\
   &=&\ \{u\in H^d(\mathscr{M};F)\mid A^tu=0, C^{\ad}\wt\rho^du=0\} \\
   &=&\ \{u\in H^d(\mathscr{M};F)\mid A^tu=0, \wt\rho^du=0\}=\ker A^t_{\mmin},
\end{eqnarray*}
where in the last line we used
\[\image C_{\frac{d}{2}}(A^t)= \{\wt\rho^du\mid u\in H^d(\mathscr{M};F), A^tu=0\} \tand \eqref{e:calderon-adjoint-boundary-condition}.\]
Now \eqref{e:ker-calderon-boundary} is done.



Note that $\ker A_{\mmin}^t$ consists of smooth sections and is finite-dimensional.
Thus we can use \eqref{e:L^2-decomposition}, i.e., the decomposition in $L^2(\mathscr{M};F)$ to get our results for $s\geq \frac{d}{2}$\/:
\begin{align*}
  H^{s-\frac{d}{2}}(\mathscr{M};F)\ =&\ L^2(\mathscr{M};F)\cap H^{s-\frac{d}{2}}(\mathscr{M};F) \\
    =&\ (\image A_C\oplus \ker A^t_{\mmin})\cap H^{s-\frac{d}{2}}(\mathscr{M};F)\\
   =&\
    \bigl( \image A_C\cap H^{s-\frac{d}{2}}(\mathscr{M};F)\bigr) \oplus \ker A^t_{\mmin}\\
     =&\ \image A_{s+\frac{d}{2},C}\oplus \ker  A^t_{\mmin},
\end{align*}
where $\mathcal{D}(A_{s+\frac{d}{2},C})=\{u\in H^{s+\frac{d}{2}}(\mathscr{M};E)\mid C\wt \rho^du=0\}$ and we have used higher regularity for well-posed boundary conditions of Remark \ref{r:regular-wellposed-boundary}d.
So $\image A_{s+\frac{d}{2},C}$ is finite-codimensional in $H^{s-\frac{d}{2}}(\mathscr{M};F)$.
Since \[\image A_{s+\frac{d}{2},C}\subset \image A_{s+\frac{d}{2}}\subset H^{s-\frac{d}{2}}(\mathscr{M};F),\]
the space $\image A_{s+\frac{d}{2}}$ is also finite-codimensional and thus closed in $H^{s-\frac{d}{2}}(\mathscr{M};F)$. So we get (\ref{e:s>d/2-L^2-decomposition}).
\qed
\end{proof}

\begin{proposition}\label{p:kernel-cont-for-s-ge-dhalf}
Let $s \geq \frac{d}{2}$\/. If $\dim A^t_{b,min}=\kappa_-$ constant for all $b\in B$, then Assumption \ref{a:continuous-family-for-s-ge-dhalf}, i.e., that the family $\bigl(A_{b,s+\frac{d}{2}}\bigr)_{b\in B}$ is continuous in the operator norm, implies that the family $\bigl(\ker A_{b,s+\frac{d}{2}}\bigr)_{b\in B}$ of closed linear subspaces is continuous in $H^{s+\tfrac d2}(\mathscr{M};E)$.
\end{proposition}

\begin{proof}
Assumption \ref{a:continuous-family-for-s-ge-dhalf} implies that the graphs $\bigl(\Graph(A_{b,s+\frac{d}{2}})\bigr)_{b\in B}$ make a continuous family of closed linear subspaces of $H^{s+\tfrac d2}(\mathscr{M};E) \x H^{s-\tfrac d2}(\mathscr{M};F)$.
Here we impose the gap topology of Definition \ref{d:gap} on the space of closed linear subspaces of the product space. Actually, the two claims are equivalent
by \cite[Theorem IV.2.23 a)]{Ka95}.

For Banach spaces $X,Y$ and any bounded linear map $Q \colon X\to Y$, we recall the elementary formulae
\[
\Graph(Q)+ X\x\{0\}\ =\ X\x\image Q \ \tand \  \Graph(Q)\cap (X\x\{0\})\ =\ \ker Q\x\{0\}.
\]
Together with Lemma \ref{l:strongdecomposition}, for $X:=H^{s+\tfrac d2}(\mathscr{M};E)$, $Y:=H^{s-\tfrac d2}(\mathscr{M};F)$ and $Q$ right-Fredholm, i.e., $\image Q$ is closed and finite-codimensional, we have
\begin{equation*}
\dim \frac{X\x Y}{\Graph(Q)+X\x\{0\}}\ =\ \dim \frac{Y}{\image Q}\
\fequal{\text{for $Q=A_{b,s+\frac{d}{2}}$}}\ \dim\ker A^t_{b,\mmin}\ = \kappa_-\/.
\end{equation*}
Now we consider the two following continuous families of closed subspaces of $X\times Y$:
\newline
 $M_b:= \Graph(A_{b,s+\frac{d}{2}})$ with $b$ running in $B$ and
the constant family
\newline
  $N_b:= H^{s+\frac{d}{2}}(\mathscr{M};E)\x\{0\}$. By Lemma \ref{l:strongdecomposition}, 
\[M_b+N_b= \Graph(A_{b,s+\frac{d}{2}}) + H^{s+\frac{d}{2}}(\mathscr{M};E)\x\{0\} = H^{s+\tfrac d2}(\mathscr{M};E) \x \image A_{b,s+\frac{d}{2}}\] are closed. Then by Lemma \ref{l:closed-continuous}b, the constance of $\kappa_-$ implies that the family
\[
M_b\cap N_b= \Graph(A_{b,s+\frac{d}{2}}) \cap (H^{s+\frac{d}{2}}(\mathscr{M};E)\x\{0\}) = \ker A_{b,s+\frac{d}{2}}\x\{0\}
\]
is continuous on $B$ and so $\bigl(\ker A_{b,s+\frac{d}{2}}\bigr)_{b\in B}$, and the proposition is proved.\qed
\end{proof}

Now we turn to the Cauchy traces $\wt \rho^d(\ker A_{b,s+\frac{d}{2}})$. Note that for $s\ge \frac{d}{2}$ the Cauchy trace operator $\wt \rho^d\colon H^{s+\frac{d}{2}}(\mathscr{M};E)\to H^s(\Si;E'^d)$ is surjective and bounded.

\begin{proposition}\label{p:Cauchy-traces-varying}
Additionally to Assumption \ref{a:continuous-family-for-s-ge-dhalf}, we assume for all $b\in B$, $\dim \ker A_{b,\mmin}=\kappa_+$ constant and $\dim \ker A^t_{b,\mmin}=\kappa_-$ constant. Then the family \[\bigl(\wt \rho^d(\ker A_{b,s+\frac{d}{2}})\bigr)_{b\in B}=\bigl(\image C_s(A_b)\bigr)_{b\in B}\] makes a continuous family of closed subspaces in $H^s(\Si;E'^d)$ for all $s\geq\frac{d}{2}$\/.
\end{proposition}	

Our proof of 
Proposition \ref{p:Cauchy-traces-varying} will use the following functional-analytic estimate.
\begin{lemma}\label{l:generalized-projection}
Let $X$, $Y$ be Banach spaces, $p\colon X\to Y$ be surjective and bounded linear. Then there exist positive constants $c$ and $\bar{c}$ such that for any closed linear subspaces $M,N$ of $X$ with $M,N\>\ker p$, we have
\[
\bar{c}\delta(M,N) \leq \delta(p(M),p(N))\leq c \delta(M,N).
\]
 where $\delta(\cdot,\cdot)$ is defined in Definition \ref{d:gap}a.
\end{lemma}
\begin{proof}
Note that $\ker p$ is a closed linear subspace of $X$.
For the quotient map $q \colon X\to X/\ker p$,
we have
$\delta(M,N)=\delta(q(M),q(N))$ (cf. \cite[Lemma A.3.1(d)]{BoZh14}).
We define the induced map $\tilde p\colon X/\ker p \to Y$ by $\tilde p(x+\ker p):=p(x)$, then $p=\tilde p\circ q$.
Since the bounded linear transformation $\tilde p\colon X/\ker p\to Y$ is bijective,
the Inverse Mapping Theorem implies that $\tilde p$ is a homeomorphism. So the lemma holds.\qed
\end{proof}

\begin{proof}[of Proposition \ref{p:Cauchy-traces-varying}]
By Proposition \ref{p:kernel-cont-for-s-ge-dhalf}, the family $\bigl(\ker A_{b,s+\frac{d}{2}}\bigr)_{b\in B}$ of closed linear subspaces is continuous in $H^{s+\tfrac d2}(\mathscr{M};E)$.
For closed subspaces $M_b:=\ker A_{b,s+\frac{d}{2}},N_b:=\ker \bigl(\wt \rho^d|_{H^{s+\frac{d}{2}}(\mathscr{M};E)}\bigr)$ of $H^{s+\frac{d}{2}}(\mathscr{M};E)$,
\[
\image C_s(A_b)=\wt \rho^d(\ker A_{b,s+\frac{d}{2}})\ 
=\ \wt \rho^d(M_b+N_b).
\]	
Moreover by Lemma \ref{l:s-Amin-A-semifredholm} and this proposition's assumption, the spaces
\[
M_b\cap N_b\ =\ \ker A_{b,s+\frac{d}{2}}\cap \ker \wt \rho^d=\ker A_{b,\mmin}
\]
are of finite constant dimension $\kappa_+$ for all $b\in B$.
Since $\image C_s(A_b)$ is closed in $H^s(\Si;E'^d)$, the subspace
\[M_b+N_b = (\wt \rho^d)^{-1}\bigl(\image C_s(A_b)\bigr)\]
is closed in $H^{s+\frac{d}{2}}(\Si;E'^d)$. So by Lemma \ref{l:closed-continuous}b, the continuous variation of $M_b=\ker A_{b,s+\frac{d}{2}}$ and the constancy of the family $N_b=\ker \bigl(\wt \rho^d|_{H^{s+\frac{d}{2}}(\mathscr{M};E)}\bigr)$ imply that the family $\left(M_b+N_b\right)_{b\in B}$ is continuous. From Lemma \ref{l:generalized-projection} we get the continuous variation of $\wt \rho^d\bigl(\ker A_{b,s+\frac{d}{2}}\bigr)=\wt\rho^d(M_b+N_b)$\/.\qed
\end{proof}

Next we will provide a non-trivial jump from the continuity of the
Cauchy data spaces to the continuity of the \Calderon\ projections.
Our arguments are based on the following observation:
Given a family of bounded projections in a Banach space, if their images and kernels are continuous in the gap topology, then this family of bounded projections is continuous in the operator norm. More precisely we have
\begin{lemma}\label{l:projector-varying1}
Let $X$ be a Banach space and $B$ be a topological space. 
Let $\left(P_b\in \Bb(X)\right)_{b\in B}$ be a family of projections, that is,
$P_b^2=P_b$ for every $b\in B$. If either

\begin{equation}\label{e:convergence-Pb-Pb0}
(1) \lim_{b\to b_0}\delta(\image P_b,\image P_{b_0})\to 0 \tand \lim_{b\to b_0}\delta(\ker P_b,\ker P_{b_0})\to 0,
\end{equation}
or \begin{equation}\label{e:convergence-Pb0-Pb}
(2) \lim_{b\to b_0}\delta(\image P_{b_0},\image P_b)\to 0 \tand \lim_{b\to b_0}\delta(\ker P_{b_0},\ker P_b)\to 0;
\end{equation}
then
\begin{equation}\label{e:convergence-operatornorm}
  \lim_{b\to b_0}\|P_b-P_{b_0}\|\to 0.
\end{equation}
\end{lemma}
\begin{proof}
We will use the quantity $\gamma(\cdot,\cdot)$ in (\ref{e:mimimal-gap}) to  get the estimate of the operator norm $\|P_b-P_{b_0}\|:=\sup_{z\in X, \|z\|=1}\|(P_b-P_{b_0})z\|$.

(1) We begin to prove (\ref{e:convergence-Pb-Pb0}) $\Rightarrow$ (\ref{e:convergence-operatornorm}).
First we recall the definition and properties of $\gamma(\cdot,\cdot)$\/.
Since $X=\image P_b\oplus \ker P_b$,
by the definition of $\gamma(\cdot,\cdot)$ in (\ref{e:mimimal-gap}), for any $x'\in\image P_b\/, y'\in \ker P_b$, we have
\begin{equation}\label{e:minimal-gap2}
  \|x'+y'\|\geq \|x'\| \gamma(\image P_b, \ker P_b) \tand  \|x'+y'\|\geq \|y'\| \gamma(\ker P_b, \image P_b);
\end{equation}
and $\gamma(\image P_b, \ker P_b)>0$, $\gamma(\ker P_b,\image P_b)>0$ (cf. \cite[Theorem IV.4.2]{Ka95}).

Then we use $\delta(\cdot,\cdot)$ and $\gamma(\cdot,\cdot)$ to give the estimate of the norm $\|P_b-P_{b_0}\|$.
Take $\delta_1:=\delta(\image P_b,\image P_{b_0})$, $\delta_2:=\delta(\ker P_b,\ker P_{b_0})$.
By the definition of $\delta(\cdot,\cdot)$ (see also \cite[IV (2.3)]{Ka95}), for any $\varepsilon>0$
 and any $z'=x'+y'$ with $x'\in \image P_b$, $y'\in \ker P_b$, we can correspondingly choose $x\in \image P_{b_0}$, $y\in \ker P_{b_0}$ such that
\begin{equation}\label{e:x'-y'-x'prime-y'prime}
  \|x'-x\|\leq (\delta_1+\varepsilon)\|x'\|,\quad \|y'-y\|\leq (\delta_2+\varepsilon)\|y'\|.
\end{equation}
So we have
\begin{align*}
 \ &\ \|(P_b-P_{b_0})z'\|
  \ =\ \|(P_b-P_{b_0})(x'+y')\|\\
  =&\ \|x'-P_{b_0}(x'+y')+P_{b_0}(x+y)-x\|\\
  =&\ \|x'-x+P_{b_0}(x'+y')-P_{b_0}(x+y)\|\\
  \leq& \ \|x'-x\|+\|P_{b_0}\|(\|x'-x\|+\|y'-y\|)\\
  \leq& \ (\|P_{b_0}\|+1)(\delta_1+\delta_2+\varepsilon)(\|x'\|+\|y'\|)\  \text{ by \eqref{e:x'-y'-x'prime-y'prime}}\\
  \leq&\ (\|P_{b_0}\|+1)(\delta_1+\delta_2+\varepsilon) (\frac{\|x'+y'\|}{\gamma(\image P_b,\ker P_b)}
  +\frac{\|x'+y'\|}{\gamma(\ker P_b,\image P_b)})\ \text{ by \eqref{e:minimal-gap2}}.
 \end{align*}
 Since $\varepsilon>0$ and $z'\in X$ are both arbitrary, we have
 \begin{equation}\label{e:estimate-operatornorm-projection}
 \|P_b-P_{b_0}\|\leq (\|P_{b_0}\|+1)(\delta_1+\delta_2)\left(\frac{1}{\gamma(\image P_b,\ker P_b)}+\frac{1}{\gamma(\ker P_b,\image P_b)}\right).
\end{equation}

Finally, we give the positive lower bound estimate of $\gamma(\image P_b,\ker P_b)$.
By \cite[Lemma 1.4]{Ne68}, if $\gamma(\image P_{b_0},\ker P_{b_0})-\delta_1\cdot
\gamma (\image P_{b_0},\ker P_{b_0})-\delta_1-\delta_2>0$, then
\[\gamma (\image P_b,\ker P_b) \geq \frac{\gamma(\image P_{b_0},\ker P_{b_0})-\delta_1\cdot
\gamma (\image P_{b_0},\ker P_{b_0})-\delta_1-\delta_2}{1+\delta_2}.\]
Together with \eqref{e:convergence-Pb-Pb0}, we have
\begin{equation}\label{e:lowerbound-gamma-imageP-kerP}
\liminf_{b\to b_0}\gamma (\image P_b,\ker P_b)\ge \gamma (\image P_{b_0},\ker P_{b_0})>0.
\end{equation}
Similarly,
\begin{equation}\label{e:lowerbound-gamma-kerP-imageP}
\liminf_{b\to b_0}\gamma (\ker P_b,\image P_b)\ge \gamma (\ker P_{b_0},\image P_{b_0})>0.
\end{equation}
Combining \eqref{e:estimate-operatornorm-projection}, \eqref{e:lowerbound-gamma-imageP-kerP} and \eqref{e:lowerbound-gamma-kerP-imageP},
 we get \eqref{e:convergence-operatornorm}.

(2) Now we are going to prove (\ref{e:convergence-Pb0-Pb}) $\Rightarrow$ (\ref{e:convergence-operatornorm}). Take $\delta_3:=\delta(\image P_{b_0},\image P_b)$, $\delta_4:=\delta(\ker P_{b_0},\image P_b)$. Similar to \eqref{e:estimate-operatornorm-projection}, we have

\begin{equation}\label{e:estimate-operatornorm-projection2}
\|P_b-P_{b_0}\| \leq (\|P_b\|+1)(\delta_3+\delta_4)
\left(\frac{1}{\gamma(\image P_{b_0},\ker P_{b_0})}+\frac{1}{\gamma(\ker P_{b_0},\image P_{b_0})}\right).
\end{equation}
Take $\alpha:=\frac{1}{\gamma(\image P_{b_0},\ker P_{b_0})}+\frac{1}{\gamma(\ker P_{b_0},\image P_{b_0})}$.
Since $\|P_b\|\leq\|P_b-P_{b_0}\|+\|P_{b_0}\|$, we have
\[\|P_b\|(1-\alpha(\delta_3+\delta_4))\leq \|P_{b_0}\|+\alpha (\delta_3+\delta_4).
\]
Together with \eqref{e:convergence-Pb0-Pb} and \eqref{e:estimate-operatornorm-projection2}, we get \eqref{e:convergence-operatornorm}.\qed
\end{proof}

By the preceding lemma and the definition of the gap, we can conclude
 \begin{corollary}\label{c:sufficient-condition-families-of-projections}
A sufficient and necessary condition for the continuity of a family of bounded projections in a fixed Banach space, parameterized by a topological space, is that their kernels and images are both continuous in the gap topology.
\end{corollary}

\begin{proof}[of Theorem \ref{t:main} for $s\geq\frac{d}{2}$]
According to Corollary \ref{c:l2-orthogonality-of-calderon-in-seeley}, we have, for $s\geq \frac{d}{2}$
 \[\ker C^{\ort}_s(A)=\wt{J}^t\left(\image C_s(A^t)\right).\] So under the assumptions of Theorem \ref{t:main}, by Proposition \ref{p:Cauchy-traces-varying}, we get that, for $s\geq \frac{d}{2}$,
 $\bigl(\image C^{\ort}_s(A_b)\bigr)_{b\in B}$ and $\bigl(\ker C^{\ort}_s(A_b)\bigr)_{b\in B}$ are both continuous in $H^s(\Si;E'^d)$. Then by Corollary \ref{c:sufficient-condition-families-of-projections}, we get that the family
 $\bigl(C^{\ort}_s(A_b)\bigr)_{b\in B}$ is continuous in the operator norm $\|\cdot\|_{s,s}$ for all $s\geq \frac{d}{2}$.
 \qed
 \end{proof}

\subsection{Proof of our main theorem for $s< \frac{d}{2}$}\label{ss:s<halfd}
Interpolation theory can be applied easily for intermediate Sobolev spaces between two given Sobolev spaces to establish an estimate for the operator norm of an intermediate operator, see \textsc{\Calderon}'s \cite{Cal64-intermediate} or \cite{LM72} by \textsc{J.-L. Lions} and \textsc{Magenes}.

We give a slimmed-down version of interpolation theory for intermediate spaces.
\begin{definition}[Interpolation property] We follow \cite[Definitions 21.4 and 21.5]{Tar07}.
	Let $\EE_0$ and $\EE_1$ be normed spaces with $\EE_1\hookrightarrow \EE_0$ continuously embedded and dense.
	\newline (a) An {\em intermediate space} between $\EE_1$ and $\EE_0$ is any normed space $\EE$ such that $\EE_1 \< \EE \< \EE_0$ (with continuous embeddings).
	\newline (b) An {\em interpolation space} between $\EE_1$ and $\EE_0$ is any intermediate space $\EE$ such that every linear mapping from $\EE_0$ into itself which is continuous from $\EE_0$ into itself and from $\EE_1$ into itself is automatically continuous from $\EE$ into itself. It is said to be of {\em exponent} $\theta$ (with $0 < \theta < 1$), if there exists a constant $c_1$ such that
	\begin{equation}\label{e:interpolation}
	\norm{A}_{\Bb(\EE,\EE)}\ \le\ c_1\, \norm{A}_{\Bb(\EE_1,\EE_1)}^{1-\theta}\,  \norm{A}_{\Bb(\EE_0,\EE_0)}^{\theta}\ \text{ for all $A\in\Bb(\EE_1,\EE_1)\cap\Bb(\EE_0,\EE_0)$}.
	\end{equation}
\newline (c)
Moreover, if $\EE_0$ and $\EE_1$ are Banach spaces, for $0 < \theta < 1$, we can define the \textit{complex interpolation space} $[\EE_1,\EE_0]_{\theta}$ in loc. cit.
\end{definition}	
\begin{remark}
The construction of the complex interpolation space uses analytic functions with values in the Banach space
$\EE_0$. Using the classical Three Lines Theorem (mainly about the maximum modulus principle), one can show that $[\EE_1,\EE_0]_{\theta}$ with a kind of quotient norm is also a Banach space (cf. \cite[Section 1.14.1]{LM72} or \cite[Section 4.2]{Taylor96}). By \cite[Lemma 21.6]{Tar07}, the interpolation property holds for $[\EE_1,\EE_0]_{\theta}$ with $c_1=1$ in \eqref{e:interpolation}\/.
\end{remark}

\begin{definition}\label{d:sobolev scale}(cf. \cite[Definition 2.5]{BrLe01})
Slightly more generally, we call a family $(H^s)_{s\in\RR}$
a \textit{scale of Hilbert spaces} if
\begin{description}
  \item[(1)] $H^s$ is a Hilbert space for each $s\in\RR$,
  \item[(2)] $H^{s'}\hookrightarrow H^s$ embeds continuously for $s\leq s'$,
  \item[(3)] if $s<t$, $0<\theta<1$, then the \textit{complex interpolation space} belongs to the scale with
   \[
   [H^t,H^s]_{\theta}=H^{(1-\theta)t+\theta s},
  \]
  \item[(4)] $H^{\infty}:=\cap_{s\in \RR_+}H^s$ is dense in $H^t$ for each $t\in\RR$,
  \item[(5)] the $H^0$-scalar product, denoted by $(\cdot,\cdot)$, restricted to $H^{\infty}$ extends to a perfect pairing between $H^s$ and $H^{-s}$, denoted by $\langle\cdot,\cdot\rangle_{s,-s}$\/, for all $s\in\RR$.
\end{description}
\end{definition}

Let $(H^s)_{s\in \RR}$ be a scale of Hilbert spaces.
\begin{definition}\label{d:0order-operator}
 A linear map $T\colon H^\infty\to H^\infty$ is called an \textit{operator of order $0$}, if it extends to a continuous linear map $T_s\colon H^s \to H^s$ for all $s\in \RR$. We denote the vector space of all operators of order $0$ by $\operatorname{Op}^0((H^s)_{s\in\RR})$\/. For $T_s\in \Bb(H^s)$, we denote its operator norm by $\|T_s\|_{s,s}$\/.
\end{definition}


\begin{lemma}\label{l:s<d-half}
Let $(H^s)_{s\in\RR}$ be a scale of Hilbert spaces  and  $T\in \operatorname{Op}^0((H^s)_{s\in\RR})$. We
assume that the continuous extension $T_0 \colon H^0\to H^0$ is self-adjoint.
Then
\begin{description}
  \item[(1)] for $t>0$, $\|T_{-t}\|_{-t,-t}=\|T_t\|_{t,t}$\,;
  \item[(2)] for $s_0<s<s_1$, $\|T_s\|_{s,s}\leq\ (\|T_{s_1}\|_{s_1,s_1})^{\tfrac{s-s_0}{s_1-s_0}}\,
      (\|T_{s_0}\|_{s_0,s_0})^{\tfrac{s_1-s}{s_1-s_0}}$\/.
\end{description}
\end{lemma}
\begin{proof}
(1) Fix any $s\in\RR$. Let $(H^s)^*$ denote the the space of bounded linear functionals on $H^s$.
The norm of $\phi\in(H^s)^*$ is given by
\[\|\phi\|:=\sup_{f\in H^s, \|f\|_s\leq 1}|\phi(f)|.\]
By Definition \ref{d:sobolev scale}(5),
$H^{-s}$ can be identified with $(H^s)^*$ by the isometric isomorphism,
\begin{eqnarray*}
  H^{-s} &\rightarrow& (H^s)^*, \\
  h &\mapsto& \langle\cdot,h\rangle_{s,-s},
\end{eqnarray*}
where isometric means that: if $\phi(f):=\langle f,h\rangle_{s,-s}$\/ for every $f\in H^s$, then $\|\phi\|=\|h\|_{-s}$.
According to the above identification,
we can define the adjoint operator of $T_s$
\begin{equation}\label{e:adjiont-of-T}
  \begin{gathered}
    (T_s)^* \colon H^{-s}\to H^{-s}\ \text{ by setting for $h\in H^{-s}$}\\
       \lla f,(T_s)^*h\rra_{s,-s}\ :=\ \lla T_sf,h\rra_{s,-s}\quad \text{for all $f\in H^s$},
   \end{gathered}
\end{equation}
then $(T_s)^*$ is also a bounded linear operator and
\begin{equation}\label{e-T-adjoint}
  \|(T_s)^*\|_{-s,-s}=\|T_s\|_{s,s}\/.
\end{equation}

For $t>0$, we claim that $(T_t)^*=T_{-t}$\/.
In fact, since $t>0$, $T_0|_{H^t}=T_t$\/, $T_{-t}|_{H^0}=T_0$\/, and for any $f\in H^t\subset H^0$, $h\in H^0\subset H^{-t}$, we have
\begin{equation}\label{e:selfadjoint-T0}
  \lla T_tf,h\rra_{t,-t}=(T_tf,h)=(T_0f,h)=(f,T_0h)=\lla f,T_0h\rra_{t,-t}=\lla f,T_{-t}h\rra_{t,-t}\/,
\end{equation}
where we have used the assumption that $T_0  \colon H^0\to H^0$ is self-adjoint.
So by \eqref{e:adjiont-of-T} and \eqref{e:selfadjoint-T0}, for any $f\in H^t,h\in H^0$, we have
\[
\lla f,(T_t)^*h\rra_{t,-t}\ =\ \lla f,T_{-t}h\rra_{t,-t}\/.
\]
This implies
\[
(T_t)^*h=T_{-t}h\ \ \  \text{for any $h\in H^0$}.
\]
Since $(T_t)^*,T_{-t}$ are
bounded linear operators on $H^{-t}$ and since $H^0$ is dense in $H^{-t}$, we get
\begin{equation}\label{e:T-adjoint-negative}
  (T_t)^*=T_{-t} \ \text{ on $H^{-t}$}.
\end{equation}
Finally by \eqref{e-T-adjoint} and \eqref{e:T-adjoint-negative}, we get $\|T_{-t}\|_{-t,-t}=\|(T_t)^*\|_{-t,-t}=\|T_t\|_{t,t}$\/.
\par
(2) Since $[H^{s_1},H^{s_0}]_{\theta}=H^{(1-\theta)s_1+\theta s_0}$, for $0<\theta<1$,
by the interpolation property for $[H^{s_1},H^{s_0}]_{\theta}$ (cf. \eqref{e:interpolation}) , we obtain, for $s_0<s<s_1$
\begin{equation*}\label{interpolation theory}
  \|T_s\|_{s,s}\\ \leq\ (\|T_{s_1}\|_{s_1,s_1})^{\tfrac{s-s_0}{s_1-s_0}}\,(\|T_{s_0}\|_{s_0,s_0})^{\tfrac{s_1-s}{s_1-s_0}}.
\end{equation*} \qed
\end{proof}

\begin{theorem}\label{t:s<d-half}
Let $B$ be a topological space and $T_b\in \operatorname{Op}^0((H^s)_{s\in\RR})$ for all $b\in B$. Assume that the extended bounded linear maps $T_{b,0} \colon H^0\to H^0$ are self-adjoint for all $b\in B$.
If $\bigl(T_{b,t}\in \Bb(H^t)\bigr)_{b\in B}$ is continuous on $B$ in the operator norm for some $t\in \RR_+$, then $\bigl(T_{b,s}\in \Bb(H^s)\bigr)_{b\in B}$ is continuous on $B$ in the operator norm for all $s\in [-t,t]$.
\end{theorem}
\begin{proof}
For any $b_1,b_2\in B$, the linear map $T_{b_1}-T_{b_2}\in \operatorname{Op}^0((H^s)_{s\in\RR})$.
According to Lemma \ref{l:s<d-half}, we have
\[
  \|T_{b_1,-t}-T_{b_2,-t}\|_{-t,-t}=\|T_{b_1,t}-T_{b_2,t}\|_{t,t}\/,
  \]
and
\[
  \|T_{b_1,s}-T_{b_2,s}\|_{s,s}\leq\  (\|T_{b_1,s_1}-T_{b_2,s_1}\|_{s_1,s_1})^{\tfrac{s-s_0}{s_1-s_0}}\,
  (\|T_{b_1,s_0}-T_{b_2,s_0}\|_{s_0,s_0})^{\tfrac{s_1-s}{s_1-s_0}}\/,
\]
where $s_0\leq s\leq s_1$ and we use the situation $s_0:=-t,s_1:=t$. So for any $s\in[-t,t]$, the continuity
of $\bigl(T_{b,s}\in \Bb(H^s)\bigr)_{b\in B}$ on $B$ in operator norm follows.\qed
\end{proof}

For the chain of Sobolev spaces over our closed manifold $\Sigma$ and $s_0 < s_1$ we set $\EE_0\ :=\  H^{s_0}(\Sigma;E'^d)$ and $\EE_1\ :=\  H^{s_1}(\Sigma;E'^d)$. 
We exploit that the Sobolev spaces are Hilbert (or Hilbertable) spaces and admit a densely defined self-adjoint positive operator $\Lambda$ in $\EE_0$\/ with domain $\mathcal{D}(\Lambda)=\EE_1$.

\begin{proposition}[Interpolation between Sobolev spaces]\label{p:interpolation}
	For each $s\in ]s_0,s_1[$ the Sobolev space $H^{s}(\Sigma;E'^d)$ is an interpolation space between
$\EE_1:=H^{s_1}(\Sigma;E'^d)$ and $\EE_0:=H^{s_0}(\Sigma;E'^d)$ of exponent
	\[
	\theta(s)\ =\ 	\frac{s_1-s}{s_1-s_0}\,.
	\]
	More precisely, we have for all $\theta\in ]0,1[$ and corresponding $s=(1-\theta)s_1+\theta s_0$\/:
	\begin{description}
		\item [(1) Identifying Sobolev spaces with interpolation spaces,]\cite[Definition 1.2.1 and Section 1.7.1]{LM72}:
		$H^{s}(\Sigma;E'^d)=\mathcal{D}(\Lambda^{1-\theta})=[\EE_1,\EE_0]_{\theta}$ with equivalent norms. The norm on $[\EE_1,\EE_0]_{\theta}$ is equivalent to the graph norm of $\Lambda^{1-\theta}$, i.e., $\bigl(\norm{u}_{\EE_0}^2 + \norm{\Lambda^{1-\theta} u}_{\EE_0}^2\bigr)^{1/2}$\/.
		\item [(2) Interpolation property of (Sobolev) norms,]\cite[Proposition 1.2.3]{LM72}: There exists a constant $c$ such that 
		$\norm{u}_{[\EE_1,\EE_0]_{\theta}}\
		\le\ c\, \norm{u}_{\EE_1}^{1-\theta}\, \norm{u}_{\EE_0}^{\theta}$\ for all $u\in \EE_1$\/.
		
		
	\end{description}
\end{proposition}
\begin{remark}
For our Hilbert spaces we have $[\EE_1,\EE_0]_{\theta}=\mathcal{D}(\Lambda^{1-\theta})$. The proof can be found in \cite[Theorem 1.14.1]{LM72} or \cite[Section 4.2]{Taylor96}. Then statements (1), (2) are immediate from the definition of the Sobolev spaces; for (2) see also \cite[Theorem 7.22]{Grubb:2009} with \textsc{Grubb}'s four-line proof in the Euclidean case based on the H\"older Inequality. 
\end{remark}

According to Proposition \ref{p:interpolation} and the facts about the chain of Sobolev spaces over a closed manifold (cf. Section \ref{ss:sobolev}), the family $\bigl( H^s(\Sigma;E'^d)\bigr)_{s\in\RR}$ satisfies
Definition \ref{d:sobolev scale}.

\begin{proof}[of Theorem \ref{t:main} for $s<\frac{d}{2}$]
We set $H^s=H^s(\Sigma;E'^d)$, $s\in\RR$, and $T_b=C^{\ort}(A_b)$\/, $b\in B$ in Theorem \ref{t:s<d-half}. By the continuity results for $s\geq \frac{d}{2}$ in Section \ref{ss:proof-for-s-ge-d-half}, we obtain our Main Theorem.\qed
\end{proof}
%
%
\section*{Appendix: Weaker conditions than Assumption (ii)}
In this Appendix, we will prove that Assumption (ii) in Theorem \ref{t:main} can be weakened a little by finer analysis above.
First, we give a kind of example about special perturbations of formally self-adjoint elliptic operator:
\begin{theorem}\label{c:appendix--A-b}
Let $A\colon \Ci(\mathscr{M};E)\to\Ci(\mathscr{M};E)$ be a formally self-adjoint elliptic operator of order $d$, i.e., $A=A^t$.
Denote by $I \colon E\to E$ the identity bundle map. Then for any $s\in \RR$, the family of $L^2$-orthogonalized Calder{\'o}n projections $\bigl(C^{\ort}_s(A-bI)\bigr)_{b\in \CC}$ is continuous at $b=0$ in the operator norm of the corresponding Sobolev space $H^s(\Si;E'^d)$.
\end{theorem}
\begin{proof}
According to Theorem \ref{t:s<d-half} and Proposition \ref{p:interpolation}, we only need to consider the case
$s\geq \frac{d}{2}$.
Then by Corollary \ref{c:l2-orthogonality-of-calderon-in-seeley}, Lemmas \ref{l:generalized-projection} and \ref{l:projector-varying1}, we only need to prove for $s\geq \frac{d}{2}$
\begin{equation}\label{e:convergence-ker+rho}
  \lim_{b\to0}\delta\bigl(\ker (A_{s+\frac{d}{2}}-bI)+\ker \wt\rho^d,\ker A_{s+\frac{d}{2}}+\ker \wt\rho^d\bigr)=0,
\end{equation}
where $\ker (A_{s+\frac{d}{2}}-bI)+\ker \wt\rho^d=(\wt \rho^d)^{-1}\bigl(\image C_s(A-bI)\bigr)$, $\ker (A_{s+\frac{d}{2}}-bI)$ and $\ker \wt\rho^d$ are all closed subspaces of $H^{s+\frac{d}{2}}(\mathscr{M};E)$.

Let $s\geq \frac{d}{2}$\/.
Since $\ker A_{s+\frac{d}{2}}\times \{0\}=\Graph(A_{s+\frac{d}{2}})\cap (H^{s+\frac{d}{2}}(\mathscr{M};E)\times \{0\})$ and $\Graph(A_{s+\frac{d}{2}})+H^{s+\frac{d}{2}}(\mathscr{M};E)\times\{0\}=H^{s+\frac{d}{2}}(\mathscr{M};E)\times \image A_{s+\frac{d}{2}}$ is closed,
by \cite[Proposition A.3.5a]{BoZh14},
we have
\[\delta\bigl(\ker (A_{s+\frac{d}{2}}-bI),\ker A_{s+\frac{d}{2}}\bigr)\leq \frac{2\delta\bigl(\Graph(A_{s+\frac{d}{2}}-b),\Graph(A_{s+\frac{d}{2}})\bigr)}{\gamma(\Graph(A_{s+\frac{d}{2}}),H^{s+\frac{d}{2}}(\mathscr{M};E) \times \{0\})}.\]
So we get
\begin{equation}\label{e:convergence-kerA-bI}
  \lim_{b\to0}\delta\bigl(\ker (A_{s+\frac{d}{2}}-bI),\ker A_{s+\frac{d}{2}}\bigr)=0.
\end{equation}
Again by \cite[Proposition A.3.5a]{BoZh14}, we have
\begin{equation}\label{e:convergence-ker-cap-ran}
  \lim_{b\to0}\delta\bigl(\ker (A_{s+\frac{d}{2}}-bI)\cap \image A_{s+\frac{3d}{2}},\ker A_{s+\frac{d}{2}}\cap \image A_{s+\frac{3d}{2}}\bigr)=0
\end{equation}
Since $C:=C^{\ort}(A)$ is a well-posed boundary condition, by Lemma \ref{l:strongdecomposition},
we have $\image A_{s+\frac{d}{2}}=\image A_{{s+ \frac{d}{2}},C}$ and then for $b\neq 0$
\begin{align*}
  \ker (A_{s+\frac{d}{2}}-bI)\subseteq \image A_{s+ \frac{d}{2}}\cap H^{s+\frac{d}{2}}(\mathscr{M};E) &=\image A_{s+ \frac{d}{2},C}\cap H^{s+\frac{d}{2}}(\mathscr{M};E) \\
  &=\image A_{s+ \frac{3d}{2},C}= \image A_{s+\frac{3d}{2}}.
\end{align*}
So for $b\neq 0$,
\begin{equation}\label{e:kerA-bI-bnot0-imageA}
\ker (A_{s+\frac{d}{2}}-bI)\subseteq \ker (A_{s+\frac{d}{2}}-bI)\cap\image A_{{s+ \frac{3d}{2}}}.
\end{equation}
By Lemma \ref{l:strongdecomposition}, we also have
\begin{equation}\label{e:imageA-cap-kerAmin-0}
  \ker A_{s+\frac{d}{2}}\cap\image A_{s+ \frac{3d}{2}}\cap \ker \wt\rho^d=\{0\}.
\end{equation}
By \cite[Lemma 1.4]{Ne68}, (\ref{e:kerA-bI-bnot0-imageA}), (\ref{e:convergence-ker-cap-ran}) and (\ref{e:imageA-cap-kerAmin-0}),
we have
\[\delta(\ker (A_{s+\frac{d}{2}}-bI)+\ker \wt\rho^d,\ker A_{s+\frac{d}{2}}+\ker \wt\rho^d)
\leq \frac{\delta(\ker (A_{s+\frac{d}{2}}-bI),\ker A_{s+\frac{d}{2}})}{\gamma(\ker (A_{s+\frac{d}{2}}-bI),\ker \wt\rho^d)},\]
and
\[\liminf_{0\neq b\to0}\gamma(\ker (A_{s+\frac{d}{2}}-bI),\ker \wt\rho^d)\geq \gamma(\ker A_{s+\frac{d}{2}}\cap\image A_{{s+ \frac{3d}{2}}},\ker \wt \rho^d)>0.\]
So together with (\ref{e:convergence-kerA-bI}), we get \eqref{e:convergence-ker+rho}.
\qed
\end{proof}

In general, we will show that Assumption (ii) in Theorem \ref{t:main}
can be weakened to Assumption (ii') in the following  Theorem \ref{t:weaker(ii)}.
Let $B,\mathscr{M},\Sigma,E,F,d$, $(A_b)_{b\in B}$ be given as in Notation \ref{n:basic-notations}.
For $s\geq \frac{d}{2}$, fix $b_0\in B$ and let
\begin{align*}
  A^{-1}_{b,s+\frac{d}{2}}(\ker A^t_{b_0,\min})&=\{u\in H^{s+\frac{d}{2}}(\mathscr{M};E)\mid A_{b,s+\frac{d}{2}}u\in \ker A^t_{b_0,\min}\}, \\
  (A^t)^{-1}_{b,s+\frac{d}{2}}(\ker A_{b_0,\min})&=\{u\in H^{s+\frac{d}{2}}(\mathscr{M};F)\mid A^t_{b,s+\frac{d}{2}}u\in \ker A_{b_0,\min}\}.
\end{align*}
Clearly,
\[\ker A_{b,s+\frac{d}{2}}\subset A^{-1}_{b,s+\frac{d}{2}}(\ker A^t_{b_0,\min}),\
\ker A^t_{b,s+\frac{d}{2}}\subset (A^t)^{-1}_{b,s+\frac{d}{2}}(\ker A_{b_0,\min}).\]
Without Assumption (ii), we still have
\begin{lemma}\label{l:gap-right-A'inverse-Amin}
Let $s \geq \frac{d}{2}$\/. Assumption \ref{a:continuous-family-for-s-ge-dhalf}, i.e., that the family $\bigl(A_{b,s+\frac{d}{2}}\bigr)_{b\in B}$ is continuous in the operator norm, implies
\begin{equation}\label{e:gap-right-A'inverse-Amin}
  \hat{\delta}(\ker A_{b_0,s+\frac{d}{2}}, A^{-1}_{b,s+\frac{d}{2}}(\ker A^t_{b_0,\min}))\to 0,
  \ \ \ \text{when $b\to b_0$.}
\end{equation}
\end{lemma}
\begin{proof}
Assumption \ref{a:continuous-family-for-s-ge-dhalf} implies that the graphs $\bigl(\Graph(A_{b,s+\frac{d}{2}})\bigr)_{b\in B}$ make a continuous family of closed linear subspaces of $H^{s+\tfrac d2}(\mathscr{M};E) \x H^{s-\tfrac d2}(\mathscr{M};F)$.
For $s\ge \frac{d}{2}$ and $b\in B$,
\begin{align*}
   &\ \{(u,A_{b,s+\frac{d}{2}}u)\in H^{s+\frac{d}{2}}(\mathscr{M};E)\times  \ker A^t_{b_0,\min}\} \\
  = &\ \Graph(A_{b,s+\frac{d}{2}})\cap (H^{s+\tfrac d2}(\mathscr{M};E)\times \ker A^t_{b_0,\min}).
\end{align*}
By Lemma \ref{l:strongdecomposition},
$\Graph(A_{b_0,s+\frac{d}{2}})\cap (H^{s+\tfrac d2}(\mathscr{M};E)\times \ker A^t_{b_0,\min})=\ker A_{b_0,s+\frac{d}{2}}\times \{0\}$.

Since the following proof holds for any $s\ge \frac{d}{2}$,
we fix an $s\ge \frac{d}{2}$ and write shorthand $A_b:=A_{b,s+\frac{d}{2}}$,
$X:=H^{s+\frac{d}{2}}(\mathscr{M};E)$, $Y:=H^{s-\frac{d}{2}}(\mathscr{M};F)$.

First, we prove that Assumption \ref{a:continuous-family-for-s-ge-dhalf} implies
\begin{equation}\label{e:gap-cap-graphA+X-kerAtmin}
  \lim_{b\to b_0}\hat{\delta}(\Graph(A_b)\cap (X\times \ker A^t_{b_0,\min}),\Graph(A_{b_0})\cap (X\times \ker A^t_{b_0,\min}))\to 0.
\end{equation}
By Lemma \ref{l:closed-continuous}a, we just need prove that Assumption \ref{a:continuous-family-for-s-ge-dhalf} implies
\begin{equation}\label{e:gap-graphA+X-kerAtmin}
  \lim_{b\to b_0}\hat{\delta}(\Graph(A_b)+ (X\times \ker A^t_{b_0,\min}),\Graph(A_{b_0})+ (X\times  \ker A^t_{b_0,\min}))\to 0.
\end{equation}
In fact, for any $b\in B$,
$\Graph(A_b)+X\times \ker A^t_{b_0,\min}=X\times (\image A_b+\ker A^t_{b_0,\min})$.
So by Lemma \ref{l:strongdecomposition},
$\Graph(A_{b_0})+X\times \ker A^t_{b_0,\min}=X\times Y$.
On one hand, the closed subspace $\Graph(A_b)+X\times \ker A^t_{b_0,\min}\subset X\times Y$.
On the other hand, by \cite[Lemma 1.4]{Ne68},
\begin{align*}
   &\  \delta(\Graph(A_{b_0})+X\times \ker A^t_{b_0,\min},\Graph(A_b)+X\times \ker A^t_{b_0,\min})\\
  \leq &\  \frac{\delta(\Graph(A_{b_0}),\Graph(A_b))}{\gamma(\Graph(A_{b_0}),X\times\ker A^t_{b_0,\min})}.
\end{align*}
So we get \eqref{e:gap-graphA+X-kerAtmin}. Thus \eqref{e:gap-cap-graphA+X-kerAtmin} holds.

Then, by the definition of the gap, Assumption \ref{a:continuous-family-for-s-ge-dhalf}
and \eqref{e:gap-cap-graphA+X-kerAtmin} imply \eqref{e:gap-right-A'inverse-Amin}.
In fact, one one hand,
\[ \delta(\ker A_{b_0},A^{-1}_b(\ker A^t_{b_0,\min}))
   \leq \delta(\Graph(A_{b_0})\cap (X\times \{0\}),\Graph(A_b)\cap (X\times \ker A^t_{b_0,\min}));
  \]
on the other hand,
\begin{eqnarray*}
&&\delta( A^{-1}_b(\ker A^t_{b_0,\min}),\ker A_{b_0})\\
&\leq & \sqrt{\|A_b\|^2+1} \cdot\delta(\Graph(A_b)\cap (X\times \ker A^t_{b_0,\min}),\Graph(A_{b_0})\cap (X\times \ker A^t_{b_0,\min})).
\end{eqnarray*}
\qed
\end{proof}

We also have the following lemma analogous to Lemma \ref{l:s-Amin-A-semifredholm}:
\begin{lemma}\label{l:s-Amin-A-inverse-Atmin}
For $s\geq \frac{d}{2}$, \[A^{-1}_{b,s+\frac{d}{2}}(\ker A^t_{b_0,\min})\cap \ker \wt \rho^d=
A^{-1}_{b,d}(\ker A^t_{b_0,\min})\cap \ker \wt \rho^d\] is finite-dimensional and consists of smooth sections.
\end{lemma}
\begin{proof}
Since $\ker A_{b,\min}=\ker A_{b,s+\frac{d}{2}}\cap \ker \wt \rho^d $ and $\ker A^t_{b_0,\min}$ are both finite-dimensional, $A^{-1}_{b,s+\frac{d}{2}}(\ker A^t_{b_0,\min})\cap \ker \wt \rho^d$ is finite-dimensional for $s\geq \frac{d}{2}$.
Since $\ker A^t_{b_0,\min}\subset C^{\infty}(\mathscr{M};F)$, by the interior regularity for elliptic operators, we have
\[A^{-1}_{b,d}(\ker A^t_{b_0,\min})\cap \ker \wt \rho^d\subset \{u\in C^{\infty}(\mathscr{M};E)\mid A_bu\in \ker A^t_{b_0,\min} \tand \wt \rho^d u=0\}.\]
Obviously, we have for $s\geq \frac{d}{2}$,
\begin{align*}
&\ \{u\in C^{\infty}(\mathscr{M};E)\mid A_bu\in \ker A^t_{b_0,\min} \tand \wt \rho^d u=0\}\\
 \subset\  &\  A^{-1}_{b,s+\frac{d}{2}}(\ker A^t_{b_0,\min})\cap \ker \wt \rho^d \\
 \subset\ &\  A^{-1}_{b,d}(\ker A^t_{b_0,\min})\cap \ker \wt \rho^d
\end{align*}
So we get the equality.
\qed
\end{proof}

Moreover, we have
\begin{lemma}\label{l:dimension-Ainverse-Atmin-cap-rho}
Let $s\geq \frac{d}{2}$, (1) $\ker A_{b,\min}\subset A^{-1}_{b,s+\frac{d}{2}}(\ker A^t_{b_0,\min})\cap \ker \wt \rho^d$;
(2) Assumption \ref{a:continuous-family-for-s-ge-dhalf} implies that
$\dim (A^{-1}_{b,s+\frac{d}{2}}(\ker A^t_{b_0,\min})\cap  \ker \wt \rho^d)\leq \dim Z_{+,0}(A_{b_0})$, when $b$ in a sufficient small neighbourhood of $b_0$ in $B$.
\end{lemma}
\begin{proof}
(1) is obvious.
(2) follows from Lemma \ref{l:gap-right-A'inverse-Amin}, \cite[Proposition A.3.5a]{BoZh14} and \cite[Corollary IV.2.6]{Ka95}.
In fact, $Z_{+,0}(A_{b_0})=\ker A_{b_0,\min}=\ker A_{b_0,s+\frac{d}{2}}\cap \ker \wt \rho^d$ and
$\lim_{b\to b_0}\delta\bigl(A^{-1}_{b,s+\frac{d}{2}}(\ker A^t_{b_0,\min})\cap \ker \wt \rho^d,\ker A_{b_0,s+\frac{d}{2}}\cap \ker \wt \rho^d\bigr)\to0$.
\qed
\end{proof}
Now, we can prove
\begin{theorem}\label{t:weaker(ii)}
Assume that
	\begin{enumerate}[(i)]
	\item  for $s\ge \frac{d}{2}$\/, the two families of bounded extensions
	\[
	\bigl(A_{b, s+\frac{d}{2}}\colon
	 H^{s+\frac{d}{2}}(\mathscr{M};E) \too H^{s-\frac{d}{2}}(\mathscr{M};F)\bigr)_{b\in B}
	 \]
and
\[\bigl(A_{b, s+\frac{d}{2}}^t\colon H^{s+\frac{d}{2}}(\mathscr{M};F) \too H^{s-\frac{d}{2}}(\mathscr{M};E)\bigr)_{b\in B}
\]
	are continuous in the respective operator norms $\norm{\cdot}_{s+\frac{d}{2},s-\frac{d}{2}}$\/, and that the family of adjusted Green's forms (of Equation \eqref{e:J-adjusted}) $\bigl(\tilde J^t_{b,s} \colon H^s(\Si;F'^d) \to H^s(\Si;E'^d)\bigr)_{b\in B}$
	is continuous in the operator norm $\norm{\cdot}_{s,s}$\/;
	\item [(ii')] $\dim (A^{-1}_{b,d}(\ker A^t_{b_0,\min})\cap \ker \wt \rho^d)=\dim Z_{+,0}(A_{b_0})$ and \\
	$\dim ((A^t)^{-1}_{b,d}(\ker A_{b_0,\min})\cap \ker \wt \rho^d)=\dim Z_{-,0}(A_{b_0})$ hold for $b$ in a neighbourhood of $b_0$ in $B$.
	\end{enumerate}
	Then for any $s\in \RR$, the family of $L^2$-orthogonalized Calder{\'o}n projections $\bigl(C^{\ort}_s(A_b)\bigr)_{b\in B}$ is continuous at $b_0$ in the operator norm of the corresponding Sobolev space $H^s(\Si;E'^d)$.
\end{theorem}

\begin{proof}
According to Theorem \ref{t:s<d-half} and Proposition \ref{p:interpolation}, we only need to prove the case
$s\geq \frac{d}{2}$.
Let $s\geq \frac{d}{2}$ in the following.
By Lemmas \ref{l:gap-right-A'inverse-Amin}, \ref{l:s-Amin-A-inverse-Atmin} and \ref{l:closed-continuous}b,
Assumption (3.1) and
$\dim (A^{-1}_{b,d}(\ker A^t_{b_0,\min})\cap \ker \wt \rho^d)=\dim Z_{+,0}(A_{b_0})$ imply
\begin{equation}\label{e:convergence-Ainverse+kerrho}
\lim_{b\to b_0}\delta(A^{-1}_{b,s+\frac{d}{2}}(\ker A^t_{b_0,\min})+\ker \wt \rho^d, \ker A_{b_0,s+\frac{d}{2}}+\ker \wt \rho^d )=0.
\end{equation}
Since $\ker A_{b,s+\frac{d}{2}}\subset A^{-1}_{b,s+\frac{d}{2}}(\ker A^t_{b_0,\min})$, \eqref{e:convergence-Ainverse+kerrho} implies
\[\lim_{b\to b_0}\delta(\ker A_{b,s+\frac{d}{2}}+\ker \wt \rho^d,\ker A_{b_0,s+\frac{d}{2}}+\ker \wt \rho^d )\to0.\]
Similarly, (i) and $\dim ((A^t)^{-1}_{b,d}(\ker A_{b_0,\min})\cap \ker \wt \rho^d)=\dim Z_{-,0}(A_{b_0})$ imply
\[\lim_{b\to b_0}\delta(\ker A^t_{b,s+\frac{d}{2}}+\ker \wt \rho^d,\ker A^t_{b_0,s+\frac{d}{2}}+\ker \wt \rho^d )\to0.\]
According to Corollary \ref{c:l2-orthogonality-of-calderon-in-seeley}, we have
\[\image C_s^{\ort}(A_b)=\wt \rho^d(\ker A_{b,s+\frac{d}{2}})\ \ \tand \ \ \ker C_s^{\ort}(A_b)= \wt J^t\wt\rho^d(\ker A^t_{b,s+\frac{d}{2}}).\]
Then applying Lemmas \ref{l:generalized-projection} and \ref{l:projector-varying1}, we get the continuity of $L^2$-orthogonalized Calder{\'o}n projections in the operator norm of the corresponding Sobolev space $H^s(\Si;E'^d)$ for $s\geq \frac{d}{2}$, then by the discussion above we get the same conclusion for all $s\in \RR$.
\qed
\end{proof}

\begin{remark}

(a) Theorem \ref{c:appendix--A-b} can be seen as a direct corollary of Theorem \ref{t:weaker(ii)}.
For $A=A^t$ and $b\in\CC$, $(A-bI)^{-1}(\ker A_{\min})=\ker (A-bI)+\ker A_{\min}$.
So when $b\to 0$ and $b\neq 0$, we have $\dim((A-bI)^{-1}(\ker A_{\min})\cap\ker \wt\rho^d)=\dim \ker A_{\min}$ and $Z_{+,0}(A-bI)=\{0\}$.
(b) By Lemma \ref{l:dimension-Ainverse-Atmin-cap-rho},
Assumption (ii) in Theorem \ref{t:main} implies Assumption (ii') in Theorem \ref{t:weaker(ii)}.
\end{remark}

\bibliography{Hamiltonian2}

\def\scr{\mathcal} \def\cprime{$'$}
\def\bysame{---\thinspace}
\renewcommand*{\bfdefault}{b}
\providecommand{\MR}{\relax\ifhmode\unskip\space\fi MR }
\providecommand{\MRhref}[2]{%
  \href{http://www.ams.org/mathscinet-getitem?mr=#1}{#2}
}
\providecommand{\href}[2]{#2}
\begin{thebibliography}{10}

\bibitem{Atiyah-Singer:1971}
\textsc{M.~F. Atiyah and I.~M. Singer}, `The index of elliptic operators.
  {IV}'. \emph{Ann. of Math. (2)} \textbf{93} (1971), 119--138.

\bibitem{Latushkin-et-al:2018}
\textsc{M.~Beck, G.~Cox, C.~Jones, Y.~Latushkin, K.~McQuighan and
  A.~Sukhtayev}, `Instability of pulses in gradient reaction-diffusion systems:
  a symplectic approach'. \emph{Philos. Trans. Roy. Soc. A} \textbf{376}/2117
  (2018), 20170187, 20.

\bibitem{BCLZ}
\textsc{B.~Boo{\ss}-Bavnbek, G.~Chen, M.~Lesch and C.~Zhu}, `Perturbation of
  sectorial projections of elliptic pseudo-differential operators'. \emph{J.
  Pseudo-Differ. Oper. Appl.} \textbf{3} (\auindex{Boo{\ss}--Bavnbek,\
  B.|bind}\auindex{Chen,\ G.|bind}\auindex{Lesch,\ M.|bind}\auindex{Zhu,\
  C.|bind}2012), 49--79. \url{arXiv:1101.0067v4 [math.SP]}.

\bibitem{BoFu98}
\textsc{B.~Boo{\ss}-Bavnbek and K.~Furutani}, `The {M}aslov index: a functional
  analytical definition and the spectral flow formula'. \emph{Tokyo J. Math.}
  \textbf{21}/1 (1998), 1--34.

\bibitem{BoLe:2009}
\textsc{B.~Boo{\ss}-Bavnbek and M.~Lesch}, `The invertible double of elliptic
  operators'. \emph{Lett. Math. Phys.} \textbf{87}/1-2 (2009), 19--46.

\bibitem{BoLeZh08}
\textsc{B.~Boo{\ss}-Bavnbek, M.~Lesch and C.~Zhu}, `The {C}alder\'on
  projection: new definition and applications'. \emph{J. Geom. Phys.}
  \textbf{59}/7 (\auindex{Boo{\ss}--Bavnbek,\ B.|bind}\auindex{Lesch,\
  M.|bind}\auindex{Zhu,\ C.|bind}2009), 784--826. \url{arXiv:0803.4160v1
  [math.DG]}.

\bibitem{BoWo93}
\textsc{B.~Boo{\ss}-Bavnbek and K.~P. Wojciechowski}, \emph{Elliptic Boundary
  Problems for {D}irac Operators}, Mathematics: Theory \& Applications.
  Birkh\"auser Boston Inc., Boston, MA, \auindex{Boo{\ss}--Bavnbek,\
  B.|bind}\auindex{Zhu,\ C.|bind}1993.

\bibitem{BoZhu:2004}
\textsc{B.~Boo{\ss}-Bavnbek and C.~Zhu}, `Weak Symplectic Functional Analysis
  and General Spectral Flow Formula'. \url{arXiv:0406139 [math.DG]}.

\bibitem{BoZh14}
\bysame, `The {M}aslov index in symplectic {B}anach spaces'. \emph{Mem. Amer.
  Math. Soc.} \textbf{252}/1201 (2018), x+118. \url{arXiv:1406.0569 [math.SG]}.

\bibitem{Bo56}
\textsc{R.~Bott}, `On the iteration of closed geodesics and the {S}turm
  intersection theory'. \emph{Comm. Pure Appl. Math.} \textbf{9} (1956),
  171--206.

\bibitem{BrLe01}
\textsc{J.~Br{\"u}ning and M.~Lesch}, `On boundary value problems for {D}irac
  type operators. {I}. {R}egularity and self-adjointness'. \emph{J. Funct.
  Anal.} \textbf{185}/1 (\auindex{Br{\"u}ning,\ J.|bind}\auindex{Lesch,\
  M.|bind}2001), 1--62.

\bibitem{Cal63}
\textsc{A.~P. Calder\'on}, `Boundary value problems for elliptic equations'.
  In: \emph{Outlines of the {J}oint Soviet-American {S}ympos. on {P}artial
  {D}ifferential {E}quations ({N}ovosibirsk, August 1963)}. Acad. Sci. USSR
  Siberian Branch, Moscow, 1963, pp.~303--304.

\bibitem{Cal64-intermediate}
\bysame, `Intermediate spaces and interpolation, the complex method'.
  \emph{Studia Math.} \textbf{24} (1964), 113--190.

\bibitem{Cal76}
\bysame, \emph{Lecture Notes on Pseudo-Differential Operators and Elliptic
  Boundary Value Problems}, Cursos de Matematica. Instituto Argentino de
  Matematica, Buenos Aires,, 1976.

\bibitem{Chernoff:1973}
\textsc{P.~R. Chernoff}, `Essential self--adjointness of powers of generators
  of hyperbolic equations'. \emph{J. Funct. Anal.} \textbf{12}/4 (1973),
  401--414.

\bibitem{Cordes:1972}
\textsc{H.~O. Cordes}, `Self--adjointness of powers of elliptic operators on
  non-compact manifolds'. \emph{Math. Ann.} \textbf{195} (1972), 257--272.

\bibitem{Frey2005On}
\textsc{C.~Frey}, \emph{On Non-local Boundary Value Problems for Elliptic
  Operators}, http://d-nb.info/1037490215/34, 2005, Inaugural-Dissertation zur
  Erlangung des Doktorgrades der Mathematisch-Naturwissenschaftlichen
  Fakult{\"a}t der Universit{\"a}t zu K{\"o}ln, Thesis advisor M. Lesch.

\bibitem{Gilkey:1995}
\textsc{P.~B. Gilkey}, \emph{Invariance Theory, the Heat Equation, and the
  {A}tiyah-{S}inger Index Theorem}, second ed., Studies in Advanced
  Mathematics. CRC Press, Boca Raton, FL, \auindex{Gilkey,\ P.|bind}1995.

\bibitem{Grubb:2009}
\textsc{G.~Grubb}, \emph{Distributions and Operators}, 1 ed., Graduate Texts in
  Mathematics 252. Springer-Verlag New York, 2009.

\bibitem{Himpel-Kirk-Lesch:2004}
\textsc{B.~Himpel, P.~Kirk and M.~Lesch}, `Calder\'{o}n projector for the
  {H}essian of the perturbed {C}hern-{S}imons function on a 3-manifold with
  boundary'. \emph{Proc. London Math. Soc. (3)} \textbf{89}/1 (2004), 241--272.

\bibitem{Ho66}
\textsc{L.~H\"{o}rmander}, `Pseudo-differential operators and non-elliptic
  boundary problems'. \emph{Ann. of Math. (2)} \textbf{83} (1966), 129--209.

\bibitem{Ka95}
\textsc{T.~Kato}, \emph{Perturbation Theory for Linear Operators}, Classics in
  Mathematics. Springer-Verlag, Berlin, 1995, Reprint of the 1980 edition.

\bibitem{LM72}
\textsc{J.-L. Lions and E.~Magenes}, \emph{Non-Homogeneous Boundary Value
  Problems and Applications. {V}ol. {I}}, Die Grundlehren der mathematischen
  Wissenschaften, Band 181. Springer-Verlag, New York-Heidelberg, 1972,
  Translated from the French by P. Kenneth.

\bibitem{Ne68}
\textsc{G.~Neubauer}, `Homotopy properties of semi-{F}redholm operators in
  {B}anach spaces'. \emph{Math. Ann.} \textbf{176} (1968), 273--301.

\bibitem{Ni95}
\textsc{L.~I. Nicolaescu}, `The Maslov index, the spectral flow, and
  decomposition of manifolds'. \emph{Duke Math. J.} \textbf{80}
  (\auindex{Nicolaescu,\ L.|bind}1995), 485--533.

\bibitem{Ni97}
\bysame, `Generalized symplectic geometries and the index of families of
  elliptic problems'. \emph{Mem. Amer. Math. Soc.} \textbf{128}/609
  (\auindex{Nicolaescu,\ L.|bind}1997), 1--80.

\bibitem{Palais-Seeley:1965}
\textsc{R.~S. Palais and R.~T. Seeley}, `Chapter XVII. Cobordism invariance of
  the analytical index'. In: \emph{Seminar on the {A}tiyah-{S}inger Index
  Theorem, R. S. Palais (ed.)}. Princeton University Press, Princeton, 1965,
  pp.~285--302.

\bibitem{See66}
\textsc{R.~T. Seeley}, `Singular integrals and boundary value problems'.
  \emph{Amer. J. Math.} \textbf{88} (1966), 781--809.

\bibitem{Seeley:1968}
\bysame, `Topics in pseudo-differential operators'. In: \emph{Pseudo-{D}iff.
  {O}perators ({C}.{I}.{M}.{E}., {S}tresa, 1968)}, C.I.M.E. Summer Schools 47.
  Edizioni Cremonese, Rome, \auindex{Seeley,\ R.T.|bind}1969, pp.~167--305,
  Reprinted in 2011 by Springer-Verlag, Berlin, Heidelberg.

\bibitem{Tar07}
\textsc{L.~Tartar}, \emph{An Introduction to {S}obolev Spaces and Interpolation
  Spaces}, Lecture Notes of the Unione Matematica Italiana, vol.~3. Springer,
  Berlin; UMI, Bologna, 2007.

\bibitem{Taylor96}
\textsc{M.~E. Taylor}, \emph{Partial Differential Equations. {I}}, Applied
  Mathematical Sciences, vol. 115. Springer-Verlag, New York, 1996, Basic
  theory.

\bibitem{Treves:1}
\textsc{F.~Tr\`eves}, \emph{Introduction to Pseudodifferential and {F}ourier
  Integral Operators. {V}ol. 1}, The University Series in Mathematics. Plenum
  Press, New York-London, 1980.

\bibitem{Zhu:2006}
\textsc{C.~Zhu}, `A generalized {M}orse index theorem'. In: \emph{Analysis,
  geometry and topology of elliptic operators}. World Sci. Publ., Hackensack,
  NJ, 2006, pp.~493--540.

\end{thebibliography}
\bibliographystyle{amsplain-jl}

%
%

\end{document}